\documentclass[a4paper,fleqn]{cas-sc}

\usepackage[numbers]{natbib}
\usepackage{graphicx}
\usepackage{color,soul}
\usepackage{array}
\usepackage{amsmath}
\usepackage{amsthm}
\usepackage{amsfonts}
\usepackage{amssymb}
\usepackage{psfrag}
\usepackage{epstopdf}
\usepackage{subcaption}
\usepackage{booktabs}
\usepackage{latexsym}
\usepackage{kantlipsum}
\usepackage{balance}
\usepackage{lscape}
\usepackage{rotating}
\usepackage{multicol}
\usepackage{multirow,bigdelim}
\usepackage{color, colortbl}
\usepackage{threeparttable}
\usepackage{float}
\usepackage{tikz}
\usetikzlibrary{calc,fit,shapes,arrows,positioning,chains,arrows.meta}
\usepackage{makeshape}
\usepackage{caption}
\usepackage{floatflt}

\newcommand{\parder}[2]{\frac{\partial #1}{\partial #2}}

\newcommand{\dparder}[2]{\dfrac{\partial #1}{\partial #2}}

\newcommand{\secder}[2]{\frac{\partial^2 #1}{\partial #2^2}}

\newcommand{\bzero}{\mathbf{0}}

\newcommand{\eval}[2][\right]{\relax\ifx#1\right\relax \left.\fi#2#1\rvert}

\newcommand{\bn}{\boldsymbol{n}}

\newcommand{\bv}{\boldsymbol{v}}

\newcommand{\bB}{\boldsymbol{B}}

\newcommand{\bD}{\boldsymbol{D}}
\newcommand{\bE}{\boldsymbol{E}}

\newcommand{\bM}{\boldsymbol{M}}
\newcommand{\bN}{\boldsymbol{N}}

\newcommand{\bV}{\boldsymbol{V}}

\newcommand{\del}{\boldsymbol{\nabla}}
\newcommand{\bcross}{\boldsymbol{\times}}

\makeatletter
\newcommand{\Rmnum}[1]{\expandafter\@slowromancap\romannumeral #1@}
\makeatother
\newcommand{\bA}{\boldsymbol{A}}

\newcommand{\bJ}{\boldsymbol{J}}
\newcommand{\bK}{\boldsymbol{K}}

\newcommand{\intomega}{\int_{\varOmega}}

\newcommand{\domega}{\,d\varOmega}

\newcommand{\bH}{\boldsymbol{H}}

\newcommand{\bphi}{\boldsymbol{\phi}}

\newcommand{\phihat}{\hat{\bphi}}

\newcommand{\bAdelta}{\bA_{\delta}}

\newcommand{\bEdelta}{\bE_{\delta}}

\newcommand{\bAhat}{\hat{\bA}}
\newcommand{\bAdeltahat}{\hat{\bA_\delta}}

\newcommand{\delphi}{\del\phi}

\newcommand{\bj}{\boldsymbol{j}}

\def\tsc#1{\csdef{#1}{\textsc{\lowercase{#1}}\xspace}}
\tsc{WGM}
\tsc{QE}
\tsc{EP}
\tsc{PMS}
\tsc{BEC}
\tsc{DE}

\begin{document}
\let\WriteBookmarks\relax
\def\floatpagepagefraction{1}
\def\textpagefraction{.001}
\shorttitle{}
\shortauthors{Durgarao Kamireddy, Saurabh Madhukar Chavan and Arup Nandy}

\title [mode = title]{Comparative Performance of Novel Nodal-to-Edge finite elements over Conventional Nodal element for Electromagnetic Analysis}
\tnotemark[1]

\tnotetext[1]{Supported by Science \& Engineering Research Board (SERB), and Department of Science \& Technology (DST), Government of India, under the project IMP/2019/000276 and VSSC, ISRO through MoU No.: ISRO:2020:MOU:NO: 480.}




\author[]{Durgarao Kamireddy}[type=editor,
                        orcid=0000-0002-8546-6306]
\cormark[1]
\ead{durga176103010@iitg.ac.in}


\address[]{Department of Mechanical Engineering, Indian Institute of Technology Guwahati, Guwahati 781039, India}

\author[]{Saurabh Madhukar Chavan}
\ead{msaurabh@iitg.ac.in}

\author[]{Arup Nandy}
\ead{arupn@iitg.ac.in}

\cortext[cor1]{Corresponding author}

\begin{abstract}
	In nodal based finite element method (FEM), degrees of freedom are associated with the nodes of the element whereas, for edge FEM, degrees of freedom are assigned to the edges of the element. Edge element is constructed based on Whitney spaces. Nodal elements impose both tangential and normal continuity of vector or scalar fields across interface boundaries. But in edge elements only tangential continuity is imposed across interface boundaries, which is consistent with electromagnetic field problems. Therefore the required continuities in the electromagnetic analysis are directly obtained with edge elements whereas in nodal elements they are attained through potential formulations. Hence, while using edge elements, field variables are directly calculated but with nodal elements, post-processing is required to obtain the field variables from the potentials. Here, we present the finite element formulations with the edge element as well as with nodal elements. Thereafter, we have demonstrated the relative performances of different nodal and edge elements through a series of examples.  All possible complexities like curved boundaries, non-convex domains, sharp corners, non-homogeneous domains have been addressed in those examples. The robustness of edge elements in predicting the singular eigen values for the domains with sharp edges and corners is evident in the analysis. A better coarse mesh accuracy has been observed for all the edge elements as compared to the respective nodal elements. Edge elements are also not susceptible to mesh distortion.
\end{abstract}

\begin{keywords}
Edge element \sep Eigen value analysis \sep Electromagnetics \sep FEM
\end{keywords}

\maketitle
\section{Introduction}
\label{sec:intro}

	In the field of computational electromagnetics, FEM has been widely used in solving various electromagnetic field problems, such as eigenvalue analysis, scattering, and radiation analysis of interior and exterior domains, etc.~(\cite{harrington2001time, Jinming-finite}). The domain of interest can be discretized with finite elements which include either nodal or edge elements to implement the FEM. However, when nodal elements are used, eigenvalue analysis shows spurious modes during eigen analysis of some specific domains~(\cite{Boffi1999, Boffi2001, Cendes1991, Nedelec1980, Reddy1994}). Vector field problems require a special type of formulations due to their special continuity requirements at material interfaces. In~\cite{Nandy2016,arupthesis,Nandy2018,Agrawal2017,Nandy2018a}, the potential formulation is used in nodal framework. But this potential formulation failed in eliminating the problem of spurious eigenvalues in sharp corner objects. In~\cite{Jog2014}, mixed finite element formulation was adopted and this was successful in the case of sharp corner objects, inhomogeneous domains in 2D. However, this method failed in the case of three dimensional curved objects. In potential formulation, field variables can not be obtained directly and post-processing is required.
\par      Whitney presented a revolutionary method to address the aforementioned limitations in the 1980s by using employed edge elements, in which degrees of freedom are assigned to the edges of the finite element rather than the nodes. These elements have been constructed using curl-conforming bases. So, these elements possess tangential continuity and normal discontinuity at material interfaces. In~\cite{Webb2005}, J P Webb mentioned various important properties of edge elements. The theoretical concept, properties, and development of edge elements were published in~\cite{Nedelec1980,Barton1987}. The construction of higher order edge elements can be done in two different ways, namely, Hierarchical and Interpolatory. Different higher order edge elements were constructed using a hierarchical approach in~\cite{Seung-CheolLee2003,Ainsworth2003,Ilic2011,Andersen1999}. In~\cite{Graglia1997,Graglia1997a,Notaros2008,Lee1991}, various higher order interpolatory elements were developed and used to analyze various field problems. In hierarchical type of edge elements, within the same discretized domain both h and p-refinements are possible, whereas in interpolatory type, only p-refinement is allowed. These vector elements were used in eigenvalue analysis for various domains in~\cite{Cendes1991,Boffi2001,Boffi2010,Ainsworth2003,Garcia-Castillo2000,Bramble2005,Kamireddy_2021,Kamireddy2020}. In~\cite{Kamireddy_2021,Kamireddy2020}, authors presented a novel conversion algorithm that converts the nodal mesh data to edge element data. In~\cite{kamireddy2021novel}, the author presented a detailed description of the conversion technique for different order edge elements.
\par     We have presented the article in the following manner. In section~\ref{math_formulation}, mathematical formulation of Maxwell's electromagnetic wave equation, variational and FEM formulation in both nodal and edge element is given. The relative performance of nodal and edge elements has been compared using benchmark examples in section~\ref{numerical_study}. The effect of mesh distortion on edge elements has also been presented in section~\ref{numerical_study}. We have concluded in section~\ref{conclusions}.

\section{Mathematical Formulation}
\label{math_formulation}
\subsection{Maxwell's Equations in Electromagnetics}
The governing differential equations for electromagnetic analysis are Maxwell's equations, given in the strong form as~\cite{griff}

\begin{subequations}
\begin{gather}
\parder{\bB}{t}+\del\bcross\bE=\bzero, \label{eqmaxwell1} \\
\del\cdot\bB=0, \label{eqmaxwell2} \\
\parder{\bD}{t}-\del\bcross\bH=-\bj, \label{eqmaxwell3} \\
\del\cdot\bD=\rho, \label{eqmaxwell4}
\end{gather}
\end{subequations}		
where the electric and magnetic fields are given by $\bf{E}$ and $\bf{H}$, the electric displacement  (electric flux) is $\bD$, the magnetic induction  (magnetic flux) is $\bB$, the charge density is $\rho$, and the current density is $\bj$.
The following constitutive relations complement the above governing equations
\begin{subequations}
\begin{align}
\bD&=\epsilon\bE, \label{eqmaxw1} \\
\bB&=\mu\bH, \label{eqmaxw2}
\end{align}
\end{subequations}
where $\mu$ and $\epsilon$ are the magnetic permeability and electric permittivity, respectively. Considering that $\epsilon$ and $\mu$ are independent of time and substituting the constitutive relations into Eqns.~\eqref{eqmaxwell1},
\eqref{eqmaxwell3} and \eqref{eqmaxwell4}, and also after eliminating $\bH$ we get
\begin{equation} \label{eqmaxwell7}
\epsilon\secder{\bE}{t}+\parder{\bj}{t}+\del\bcross\left (\frac{1}{\mu}\del\bcross\bE\right)=\bzero.
\end{equation}

We get the compatibility condition from Eqns.~\eqref{eqmaxwell3}, \eqref{eqmaxw1} and \eqref{eqmaxwell4} as
\begin{equation*}
\parder{\rho}{t}+\del\cdot\bj=\bzero.
\end{equation*}

We have the boundary condition as on conducting boundary, $\bE\bcross\bn=\bzero$. Also, across the material interface, both $\bE\bcross\bn$ and $\bH\bcross\bn$ are continuous as there is no impressed surface currents. 

Introducing the relative permeability and relative permittivity $\mu_r:=\mu/\mu_0$ and $\epsilon_r:=\epsilon/\epsilon_0$, where $\mu_0$ and $\epsilon_0$ are the permeability and permittivity for the vacuum, Eqn.~\eqref{eqmaxwell7} can be written in the frequency domain as
\begin{equation} \label{eqmaxwell10}
\del\bcross\left (\frac{1}{\mu_r}\del\bcross\bE\right)-k_0^2\epsilon_r\bE=-i\omega\mu_0 \bj,
\end{equation}
where $k_0=\omega/c$ is the wave number, $c=1/\sqrt{\epsilon_0\mu_0}$, and $i=\sqrt{-1}$. Considering $\bj$ to be zero, the Eqn.~\ref{eqmaxwell10} becomes 
\begin{equation} \label{eqmaxwell11}
\del\bcross\left (\frac{1}{\mu_r}\del\bcross\bE\right)=k_0^2\epsilon_r\bE.
\end{equation}
which governs the eigenvalue problem.

\subsection{Variational statement in the nodal framework} \label{subsec_poten}

For homogeneous domains, we get wrong multiplicities of eigenvalues and for inhomogeneous domains, we get spurious values using regularized formulations in the nodal framework~(\cite{Paulsen1991}). It is partially overcome in the regularized potential formulation~(\cite{Transactions1991}), which is fairly robust in all convex domains, homogeneous or inhomogeneous. Only, it fails to find the singular eigenvalues for a non-convex domain, owing to the penalty term.

Replacing $\bE$ by $\bA+\delphi$ in Eqn.~\eqref{eqmaxwell11}, we have
\begin{equation} \label{eqmax1}
\del\bcross\left (\frac{1}{\mu_r}\del\bcross\bA\right) = k_0^2\epsilon_r\bA+k_0^2\epsilon_r\delphi,
\end{equation}
where $\phi$ and $\bA$ are scalar and vector potentials.
We get the required variational statement after adding the penalty term~(\cite{Transactions1991}) as

\begin{align}
&\intomega \frac{1}{\mu_r} (\del\bcross\bAdelta)\cdot (\del\bcross\bA)\domega + \intomega\frac{1}{\epsilon_r\mu_r} (\del\cdot\bAdelta)[\del\cdot(\epsilon_r\bA)]\domega = \notag \\
&\qquad k_0^2\intomega\epsilon_r\bAdelta\cdot\bA\domega+k_0^2\intomega\epsilon_r\bAdelta\cdot\delphi.  \label{eqmax8}
\end{align}
Multiplying Eqn.~\eqref{eqmaxwell4} by the variation $\phi_\delta$, replacing $\bE$ by $\bA+\delphi$, and considering no charge for eigen analysis we have
\begin{equation}
k_0^2\intomega \delphi_\delta\cdot (\epsilon_r\bA)\domega + k_0^2\intomega \delphi_\delta\cdot (\epsilon_r\delphi)\domega = \bzero \label{eqmax9}
\end{equation}
 On parts of the boundary where $\bE\bcross\bn=\bzero$, we specify $\phi$ and $\bA\bcross\bn$ to be zero.

\subsection{FEM Formulation in nodal framework}
The vector potential $\bA$ and its variation $\bAdelta$ are discretized as
\begin{align*}
\bA&=\bN\bAhat, \\
\bAdelta&=\bN\bAdeltahat,
\end{align*}
leading to
\begin{align*}
\del\bcross\bA&=\bB\bAhat, & \del\cdot\bA&=\bB_p\bAhat, \\
\del\bcross\bAdelta&=\bB\bAdeltahat, & \del\cdot\bAdelta&=\bB_p\bAdeltahat.
\end{align*}
Similarly, for the scalar potential $\phi$, we have
\begin{align*}
\phi&=\bN_{\phi}\phihat, & \delphi&=\bB_\phi\phihat, \\
\phi_\delta&=\bN_{\phi}\phihat_{\delta}, & \delphi_\delta&=\bB_\phi\phihat_{\delta}.
\end{align*}
where

\begin{equation*} 
\begin{aligned} 
\bN&=\begin{bmatrix} 
 N_1 & 0   & 0   &  N_2 & 0   & 0  &\ldots \\[3mm]
 0   & N_1 & 0   &  0   & N_2 & 0  &\ldots \\[3mm]
0   & 0   & N_1 &  0   & 0   & N_2 &\ldots \\[3mm] 
\end{bmatrix},
\end{aligned}
\end{equation*}
\begin{equation*} 
\begin{aligned} 
\bN_\phi=\begin{bmatrix} 
 N_1 & N_2 & N_3 &  \ldots 
 \end{bmatrix}, 
 \end{aligned}
\end{equation*}
\begin{equation*} 
\begin{aligned} 
\bB&=\begin{bmatrix} 
0 & -\dparder{N_1}{z} & \dparder{N_1}{y} &  0 & -\dparder{N_2}{z} & -\dparder{N_2}{y} &\ldots \\[3mm]
\dparder{N_1}{z} & 0 & -\dparder{N_1}{x} &  \dparder{N_2}{z} & 0 & -\dparder{N_2}{x} & \ldots \\[3mm]
-\dparder{N_1}{y} & \dparder{N_1}{x} & 0 & -\dparder{N_2}{y} & \dparder{N_2}{x} & 0 & \ldots \\[3mm] \end{bmatrix}, \\ 
\bB_p&=\begin{bmatrix} 
 \dparder{N_1}{x} & \dparder{N_1}{y} &  \dparder{N_1}{z} & \dparder{N_2}{x} & \dparder{N_2}{y} &  \dparder{N_2}{z} & \ldots \\[3mm]
 \end{bmatrix},
  \end{aligned}
\end{equation*}
  \begin{equation*} 
\quad\bB_\phi=\begin{bmatrix} 
 \dparder{N_1}{x} & \dparder{N_2}{x} &  \dparder{N_3}{x} & \ldots \\[3mm]
 \dparder{N_1}{y} & \dparder{N_2}{y} &  \dparder{N_3}{y} & \ldots \\[3mm]
\dparder{N_1}{z} & \dparder{N_2}{z} &  \dparder{N_3}{z} & \ldots \\[3mm] \end{bmatrix}. 
\end{equation*}

The discretized forms of Eqns.~\eqref{eqmax8} and~\eqref{eqmax9} are given by
\begin{equation}
\begin{bmatrix} \bK_{AA} & \bzero \\ \bzero & \bzero \end{bmatrix}
\begin{bmatrix} \bAhat \\ \phihat \end{bmatrix}=k_0^2
\begin{bmatrix} \bM_{AA} & \bM_{A\phi} \\ \bM_{\phi A} & \bM_{\phi\phi} \end{bmatrix}
\begin{bmatrix} \bAhat \\ \phihat \end{bmatrix},
\end{equation}
where

\noindent\begin{subequations}
   \begin{align}
     \bK_{AA} = \intomega \frac{1}{\mu_r}\left[\bB^T\bB+\bB_p^T\bB_p\right]\domega,  
    \end{align}
    \begin{align}
      \label{a}
      \bM_{AA} &= \intomega \epsilon_r\bN^T\bN\domega, \\
      \label{c}
      \bM_{A\phi} &= \intomega \epsilon_r\bN^T\bB_\phi\domega,
    \end{align}
    \begin{align}
     \bM_{\phi A} &= \intomega \epsilon_r\bB_\phi^T\bN\domega,  \label{eqmax13} \\
     \bM_{\phi\phi} &= \intomega \epsilon_r\bB_\phi^T\bB_\phi\domega. \label{eqmax14}
    \end{align}
\end{subequations}



\subsection{FEM formulation in edge element framework}
The variational formulation of Eqn.~\eqref{eqmaxwell11} can be derived as
\begin{equation}
\intomega \frac{1}{\mu_r} (\del\bcross\bEdelta)\cdot (\del\bcross\bE)\domega=k_0^2\intomega\epsilon_r\bEdelta\cdot\bE\domega. \label{eqmaxwell17}
\end{equation}
where the boundary conditions $\bE\bcross\bn$ and $\bH\bcross\bn$ are specified on the surfaces $\varGamma_e$ and $\varGamma_h$ of the domain respectively. For the conducting boundary, we have $\bE\bcross\bn = \bzero$ and for eigen analysis, we have $\bH = \bzero$.

We discretize the fields and their variations in Eqn.\ref{eqmaxwell17} as
\begin{align*}
\bE=\bV\hat{\bE},  &&&&   \bEdelta=\bV\hat{\bE_\delta}, \\
\del\bcross\bE=\bB\hat{\bE},    &&&&  \del\bcross\bEdelta=\bB\hat{\bE_\delta}, \\
\end{align*}
where $\hat{\bE}$ is the values of $\bE$ at different edges, $\hat{\bE_\delta}$ denote the respective variation of $\hat{\bE}$. Edge shape functions matrix, $\bV$ and $\bB$-matrix are given as
\begin{align*}
\bV &= \begin{bmatrix}
v_{1x} & v_{2x} &  \hdots \\
v_{1y} & v_{2y} &  \hdots \\
v_{1z} & v_{2z} &  \hdots \\
\end{bmatrix}, 
\quad \bB =
\begin{bmatrix}
\parder{v_{1z}}{y}-\parder{v_{1y}}{z} & \parder{v_{2z}}{y}-\parder{v_{2y}}{z} & \hdots  \\[2mm]
\parder{v_{1x}}{z}-\parder{v_{1z}}{x} & \parder{v_{2x}}{z}-\parder{v_{2z}}{x} & \hdots &  \\[2mm]
\parder{v_{1y}}{x}-\parder{v_{1x}}{y} & \parder{v_{2y}}{x}-\parder{v_{2x}}{y} & \hdots &  \\[2mm] 
\end{bmatrix}. 
\end{align*}
 where $v_x$, $v_y$ and $v_z$ are the $x$, $y$, $z$ components of each edge shape function.

After substituting the above discretizations into Eqn.\ref{eqmaxwell17} and using the arbitrariness of variations we get

\begin{equation}
\bK\hat{\bE}=k_0^2\bM\hat{\bE},
\end{equation}
where,

\noindent\begin{subequations}
\begin{minipage}{0.5\textwidth}
\begin{align}
\bK = \intomega\frac{1}{\mu_r}\bB^T\bB\domega,
\end{align}
\end{minipage}
\begin{minipage}{0.5\textwidth}
\begin{align}
\bM = \intomega\epsilon_r\bV^T\bV\domega
\end{align}
\end{minipage}
\end{subequations}
\subsection{Calculation of $\del\xi$, $\del\eta$}
$\del\xi$ and $\del\eta$ appear in different shape functions of various edge elements, hence we have to understand how they are found from the components of inverse Jacobian. We know that
\begin{equation}
\begin{Bmatrix}
\parder{f}{\xi} \\[2mm]
\parder{f}{\eta}  
\end{Bmatrix} = \bJ \begin{Bmatrix}
\parder{f}{x} \\[2mm]
\parder{f}{y}
\end{Bmatrix} \implies 
\begin{Bmatrix}
\parder{f}{x} \\[2mm]
\parder{f}{y}
\end{Bmatrix}=\bJ^{-1}\begin{Bmatrix}
\parder{f}{\xi} \\[2mm]
\parder{f}{\eta}  
\end{Bmatrix},  \label{jacobian}  
\end{equation}
where Jacobian $\bJ$ can be written as
\begin{equation*}
\bJ=\begin{bmatrix}
J_{11} & J_{12} \\[2mm] J_{21}  & J_{22} 
\end{bmatrix}=
\begin{bmatrix}
\parder{x}{\xi} & 
\parder{y}{\xi}  \\[2mm]
\parder{x}{\eta} &
\parder{y}{\eta} 
\end{bmatrix} 
\text{and assume } \bJ^{-1}=\mathbf{\Gamma}=\begin{bmatrix}
\Gamma_{11} & \Gamma_{12} \\[2mm] \Gamma_{21}  & \Gamma_{22} 
\end{bmatrix}.
\end{equation*}
 In two dimension, natural coordinates are $\xi$ and $\eta$. For $f = \xi$ in Eqn.~\ref{jacobian} we replace
$\parder{f}{\xi}=1$, $\parder{f}{\eta}=0$.
\begin{align*} \therefore
\begin{Bmatrix}
\parder{\xi}{x} \\[2mm]
\parder{\xi}{y}
\end{Bmatrix}=\begin{bmatrix}
\Gamma_{11} & \Gamma_{12} \\[2mm] \Gamma_{21}  & \Gamma_{22}
\end{bmatrix}\begin{Bmatrix}
1 \\[2mm]
0 \end{Bmatrix}= \begin{Bmatrix}
\Gamma_{11} \\[2mm]
\Gamma_{21}
\end{Bmatrix}.
\end{align*}
Similarly, for $f=\eta$ we get
\begin{align*}
\begin{Bmatrix}
\parder{\eta}{x} \\[2mm]
\parder{\eta}{y}
\end{Bmatrix}=\begin{Bmatrix}
\Gamma_{12} \\[2mm]
\Gamma_{22}
\end{Bmatrix}.
\end{align*}

\subsection{Different edge elements}
Different nodal elements used in this work are well known in the literature. We denote Four node quadrilateral elements by Q4, Nine node quadrilateral elements by Q9, and Six node triangular elements by T6. We are presenting different edge elements in this section. These elements include 4-edge quadrilateral, 12-edge quadrilateral and 8-edge triangular elements.

\begin{figure}[pos=ht]
\centering
\begin{subfigure}{0.3\textwidth}
\centering
\includegraphics[width=0.7\textwidth]{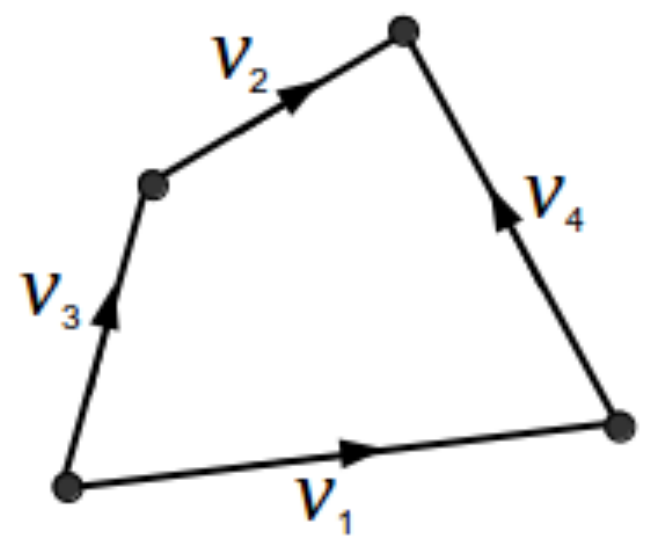}
\caption{Four edge quadrilateral element}
\label{4_edge_quad}
\end{subfigure}%
\begin{subfigure}{0.3\textwidth}
\centering
\includegraphics[width=0.7\textwidth]{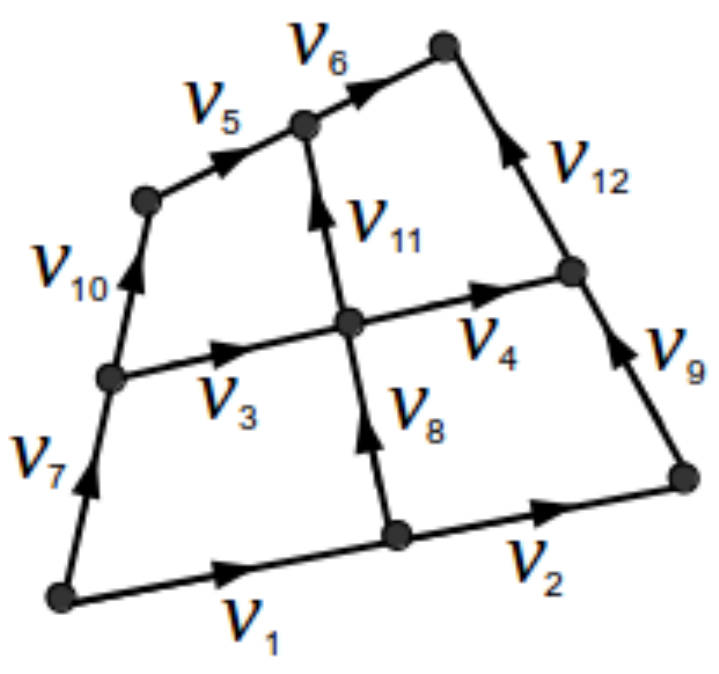}
\caption{Twelve edge quadrilateral element}
\label{12_edge_quad}
\end{subfigure}
\begin{subfigure}{0.3\textwidth}
\centering
\includegraphics[width=0.7\textwidth]{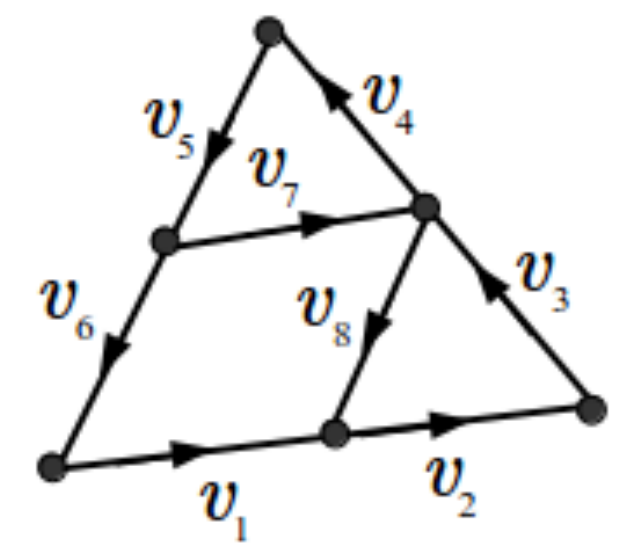}
\caption{Twelve edge quadrilateral element}
\label{8_edge_tri}
\end{subfigure}
\caption{Different edge elements}
\label{edge_elements}%
\end{figure}

\subsubsection{Four edge quadrilateral element}
Fig.~\ref{4_edge_quad} shows the quadrilateral edge element with four edges. The arrow directions are representing positive convention directions along those edges. We denote this element by EQ4. Four edge shape functions $\bv_1$, $\bv_2$, $\bv_3$, and $\bv_4$ are given as
\begin{align*}
\begin{Bmatrix}
\bv_1 \\[2mm] \bv_2 \\[2mm] \bv_3 \\[2mm] \bv_4
\end{Bmatrix}=
\begin{Bmatrix}
\frac{l_1}{4}(1-\eta)\del\xi \\[2mm]
\frac{l_2}{4}(1+\eta)\del\xi \\[2mm]
\frac{l_3}{4}(1-\xi)\del\eta \\[2mm]
\frac{l_4}{4}(1+\xi)\del\eta \\
\end{Bmatrix}
\end{align*}
where $l_1$, $l_2$, $l_3$, and $l_4$ are lengths of edges 1, 2, 3, and 4 respectively.

\subsubsection{Twelve edge quadrilateral element:}


Fig.~\ref{12_edge_quad} shows the higher order quadrilateral element with twelve edges. This element is denoted by EQ12. $\bv_1$, $\bv_2$, ..., $\bv_{11}$, and $\bv_{12}$ are the edge shape functions of edges 1 to 12 respectively. These shape functions are given as
\begin{align*}
\begin{Bmatrix}
\bv_1 \\[2mm] \bv_2 \\[2mm] \bv_3 \\[2mm] \bv_4 \\[2mm] \bv_5 \\[2mm] \bv_6 \\[2mm] \bv_7 \\[2mm] \bv_8 \\[2mm] \bv_9 \\[2mm] \bv_{10} \\[2mm] \bv_{11} \\[2mm] \bv_{12} 
\end{Bmatrix}=
\begin{Bmatrix}
\frac{-l_1}{2}\eta(\eta-1)(\xi-0.5)\del\xi \\[2mm]
\frac{l_2}{2}\eta(\eta-1)(\xi+0.5)\del\xi \\[2mm]
l_3(\eta^2-1)(\xi-0.5)\del\xi \\[2mm]
-l_4(\eta^2-1)(\xi+0.5)\del\xi \\[2mm]
\frac{-l_5}{2}\eta(\eta+1)(\xi-0.5)\del\xi \\[2mm]
\frac{l_6}{2}\eta(\eta+1)(\xi+0.5)\del\xi \\[2mm]
\frac{-l_7}{2}\xi(\xi-1)(\eta-0.5)\del\eta \\[2mm]
l_8(\xi^2-1)(\eta-0.5)\del\eta \\[2mm]
\frac{-l_9}{2}\xi(\xi+1)(\eta-0.5)\del\eta \\[2mm]
\frac{l_{10}}{2}\xi(\xi-1)(\eta+0.5)\del\eta \\[2mm]
-l_{11}(\xi^2-1)(\eta+0.5)\del\eta \\[2mm]
\frac{l_{12}}{2}\xi(\xi+1)(\eta+0.5)\del\eta
\end{Bmatrix}
\end{align*}
where $l_1$, $l_2$, ..., $l_{11}$, and $l_{12}$ are lengths of edges 1, 2, ..., 11, and 12 respectively.

\subsubsection{Eight edge triangular element}

Fig.~\ref{8_edge_tri} shows the triangular edge element with eight edges. ET8 is used to denote this edge element in this work.
 $\bv_1$, $\bv_2$, ..., $\bv_7$, and $\bv_8$ are the edge shape functions of edges 1, 2, ..., 7, and 8 edges respectively. These edge shape functions can be given as \newline
\begin{align*}
\begin{Bmatrix}
\bv_1 \\[2mm] \bv_2 \\[2mm] \bv_3 \\[2mm] \bv_4 \\[2mm] \bv_5 \\[2mm] \bv_6 \\[2mm] \bv_7 \\[2mm] \bv_8
\end{Bmatrix}=
\begin{Bmatrix}
l_1(4\xi-1)(\xi\del\eta-\eta\del\xi) \\[2mm]
l_2(4\eta-1)(\xi\del\eta-\eta\del\xi) \\[2mm]
l_3(4\eta-1)(\eta\del\alpha-\alpha\del\eta) \\[2mm]
l_4(4\alpha-1)(\eta\del\alpha-\alpha\del\eta) \\[2mm]
l_5(4\alpha-1)(\alpha\del\xi-\xi\del\alpha) \\[2mm]
l_6(4\xi-1)(\alpha\del\xi-\xi\del\alpha) \\[2mm]
4l_7\eta(\alpha\del\xi-\xi\del\alpha) \\[2mm]
4l_8\xi(\eta\del\alpha-\alpha\del\eta)
\end{Bmatrix}
\end{align*}
where $l_1$, $l_2$, ... $l_{7}$, and $l_{8}$ are lengths of edges 1, 2, ..., 7, and 8 respectively.

Partial derivatives of edge shape functions $\bv_{1}, \bv_{2} $ etc. with respect to $\xi$ and $\eta$ i.e., $\parder{v_{1x}}{\xi}, \parder{v_{1x}}{\eta}, ...  $ can be found using Mathematica~(\cite{mathematica}) for EQ4 and EQ12 whereas for ET8 they can be obtained by Finite difference method (FDM).
\begin{equation}
\parder{v_1}{\xi} = \lim_{\Delta\xi \to 0} \frac{v_{1}(\xi+\Delta\xi,\eta)-v_{1}(\xi,\eta)}{\Delta\xi} 
\end{equation}
\begin{equation} 
\parder{v_1}{\eta} = \lim_{\Delta\eta \to 0} \frac{v_{1}(\xi,\eta+\Delta\eta)-v_{1}(\xi,\eta)}{\Delta\eta}  
\end{equation}

\section{Numerical Examples}
\label{numerical_study}
\subsection{Comparative analysis of mesh convergence study between nodal elements and edge elements} 

We have done some comparative performance study among three nodal elements Q4, Q9, and T6 and three edge elements EQ4, EQ12, and ET8 through several standard benchmark examples. For all the problems discussed in this section, we have assumed $\epsilon_r=\mu_r=1.0$. In the following examples, we have compared edge and nodal elements with respect to a certain error percentage with analytical benchmark values. In that bar-diagram comparison, each bar represents the minimum no. of required FDOF for that particular element to attain an error percentage less than the predefined scale.

There are black crosses in some bar diagrams which signify that it is not possible to generate less than 10$\%$ error with that element with available computational resources. We have mentioned best possible result for those cases. Such situations mostly occur with nodal element which further establish the better performance of edge elements.


\subsubsection{Square domain }     \label{squaresection}
A square domain with a side length of $\pi$ is considered. The square domain's sides/boundaries are all perfectly conducting. Analysis data of different nodal and edge elements, like no. of free degree of freedom  (FDOF) i.e. total no. of equations, is presented in Table~\ref{tabsquare1}.

\begin{table*}[pos=h!]
\caption{Analysis data of different nodal and edge elements for the square domain problem.} \label{tabsquare1}
\begin{center}
\begin{tabular}{llllll}\hline
\multicolumn{3}{l}{{Nodal element}} &\multicolumn{3}{l}{{Edge element}} \\ \cmidrule(rr){1-3}\cmidrule(){4-6} 
\multicolumn{1}{l}{{Element}} & \multicolumn{1}{l}{{No. of}} & \multicolumn{1}{l}{{No. of}} & \multicolumn{1}{l}{{Element}} & \multicolumn{1}{l}{{No. of}} & \multicolumn{1}{l}{{No. of}} \\
\multicolumn{1}{l}{{Type}} & \multicolumn{1}{l}{{elements}} & \multicolumn{1}{l}{{FDOF/Equations}} & \multicolumn{1}{l}{{Type}} & \multicolumn{1}{l}{{elements}} & \multicolumn{1}{l}{{FDOF/Equations}}\\\cmidrule(rr){1-3} \cmidrule{4-6} 
Q4  & 64 & 231 & EQ4  & 196 & 364 \\
T6  & 98 & 663 & ET8  & 112 & 530 \\
Q9  & 64 & 855 & EQ12 & 100 & 760 \\  \hline
\end{tabular}
\end{center}
\end{table*}

\begin{table*}[pos=h!]
\begin{center}
\caption{$k_0^2$ on the square domain for different elements.} \label{tablesquare2}
\begin{tabular}{lllllll}\hline
\multicolumn{1}{l}{Analytical} & \multicolumn{3}{l}{Nodal element} &\multicolumn{3}{l}{Edge element}  \\ \hline
Benchmark &     Q4     &     Q9      &      T6      &       EQ4     &       EQ12     &      ET8    \\ \hline
1         &  1.012916  &   0.999843  &    0.999919  &    1.004203   &   1.000013  &   0.999990   \\  
1         &  1.012916  &   0.999843  &    0.999919  &    1.004203   &   1.000013  &   1.000015   \\  
2         &  2.025832  &   1.999264  &    2.000128  &    2.008407   &   2.000027  &   2.000150   \\  
4         &  4.209548  &   3.998554  &    4.000663  &    4.067583   &   4.000850  &   3.999844   \\  
4         &  4.209548  &   3.998554  &    4.001671  &    4.067583   &   4.000850  &   4.000469   \\  
5         &  5.222465  &   4.997901  &    5.003591  &    5.071786   &   5.000863  &   5.000285   \\  
5         &  5.222465  &   4.997901  &    5.007171  &    5.071786   &   5.000863  &   5.000285   \\  
8         &  8.419096  &   7.996779  &    8.033095  &    8.135166   &   8.001700  &   8.008889   \\  
9         &  10.080293 &   9.007072  &    9.024805  &    9.344778   &   9.009435  &   8.998040   \\  
9         &  10.080293  &  9.007072  &    9.024956  &    9.344778   &   9.009435  &   9.005165   \\  
10        &  11.093208  &  10.002567  &   10.038396  &   10.348982  &   10.009449  &  10.003717  \\  
10        &  11.093208  &  10.012254  &   10.048274  &   10.348982  &   10.009449  &  10.011721  \\  
13        &  14.289839  &  13.015767  &   13.088936  &   13.412361  &   13.010285  &  13.014842  \\
13        &  14.289839  &  13.015767  &   13.179304  &   13.412361  &   13.010285  &  13.048523  \\
16        &  19.453669  &  16.064853  &   16.145861  &   17.100354  &   16.051302  &  15.989497  \\
16        &  19.453669  &  16.076531  &   16.158926  &   17.100354  &   16.051302  &  16.025185  \\ \hline
\multicolumn{4}{l}{Number of computed zeros}  &  \multicolumn{3}{l}{}       \\  \midrule
-   &       77     &     245      &   221     &    45    &   121  &  65   \\ \bottomrule
\end{tabular}
\end{center}
\end{table*}

\begin{figure}[pos=h!]
\centering
\begin{subfigure}{0.49\textwidth}
\centering
\includegraphics[width=0.7\textwidth]{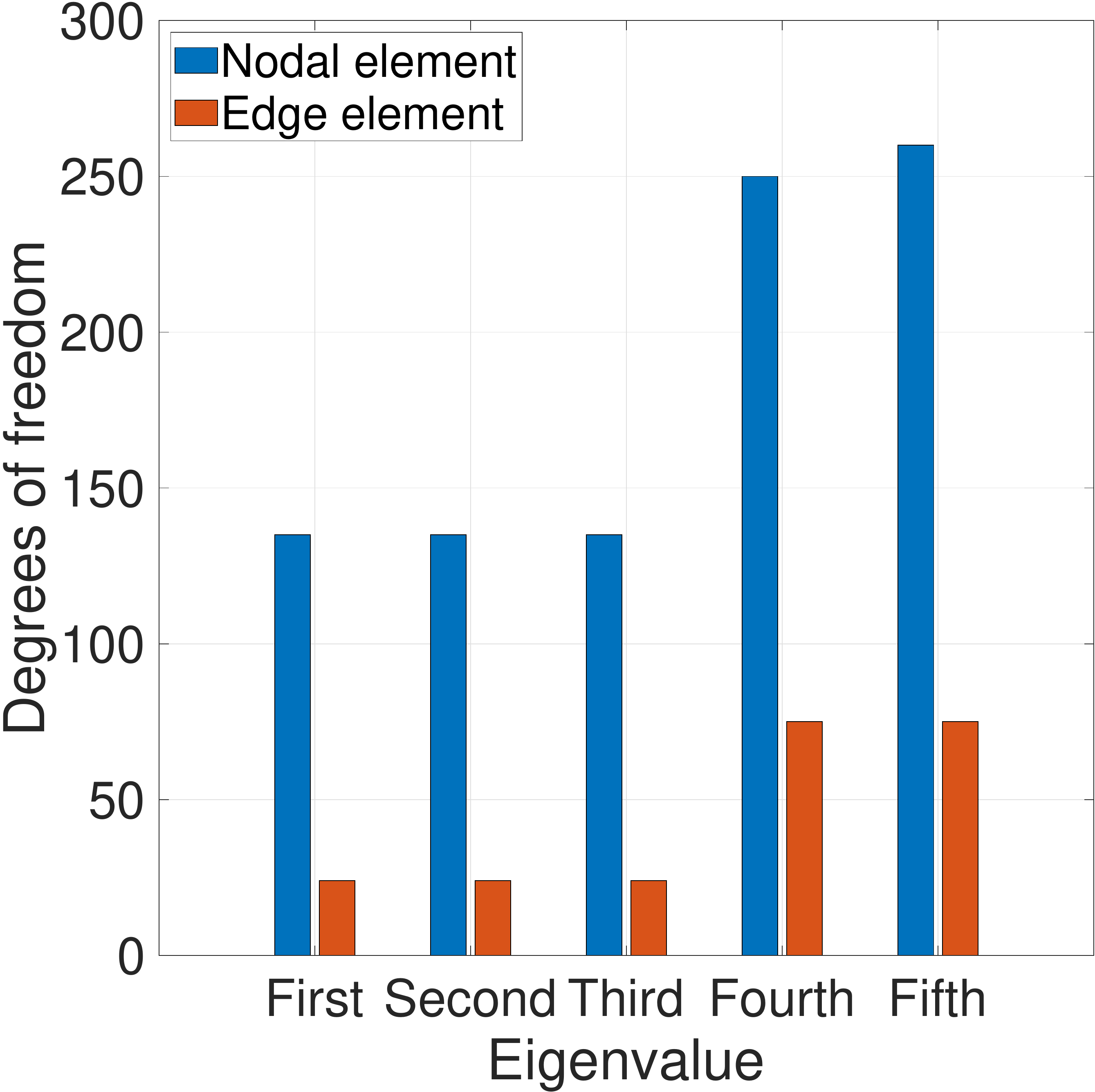}        
\caption{Minimum number of FDOF required for less than 1$\%$ error (Q9 vs EQ12).}
\label{Sa}
\end{subfigure}%
\hfill
\begin{subfigure}{0.49\textwidth}
\centering
\includegraphics[width=0.7\textwidth]{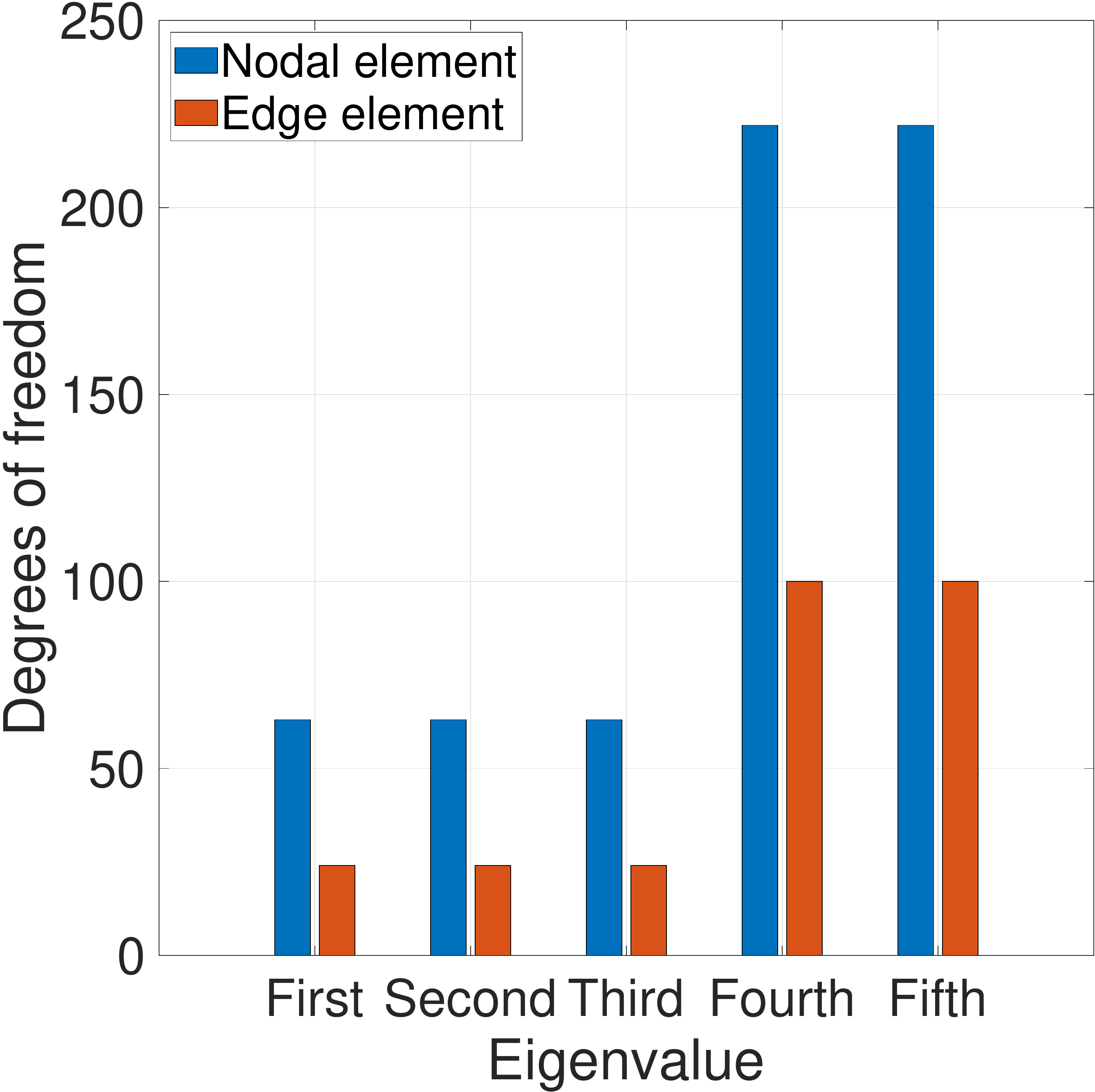}
\caption{Minimum number of FDOF required for less than 6$\%$ error (Q4 vs EQ4).}
\label{Sa}
\end{subfigure}
\begin{subfigure}{0.49\textwidth}
\centering
\includegraphics[width=0.7\textwidth]{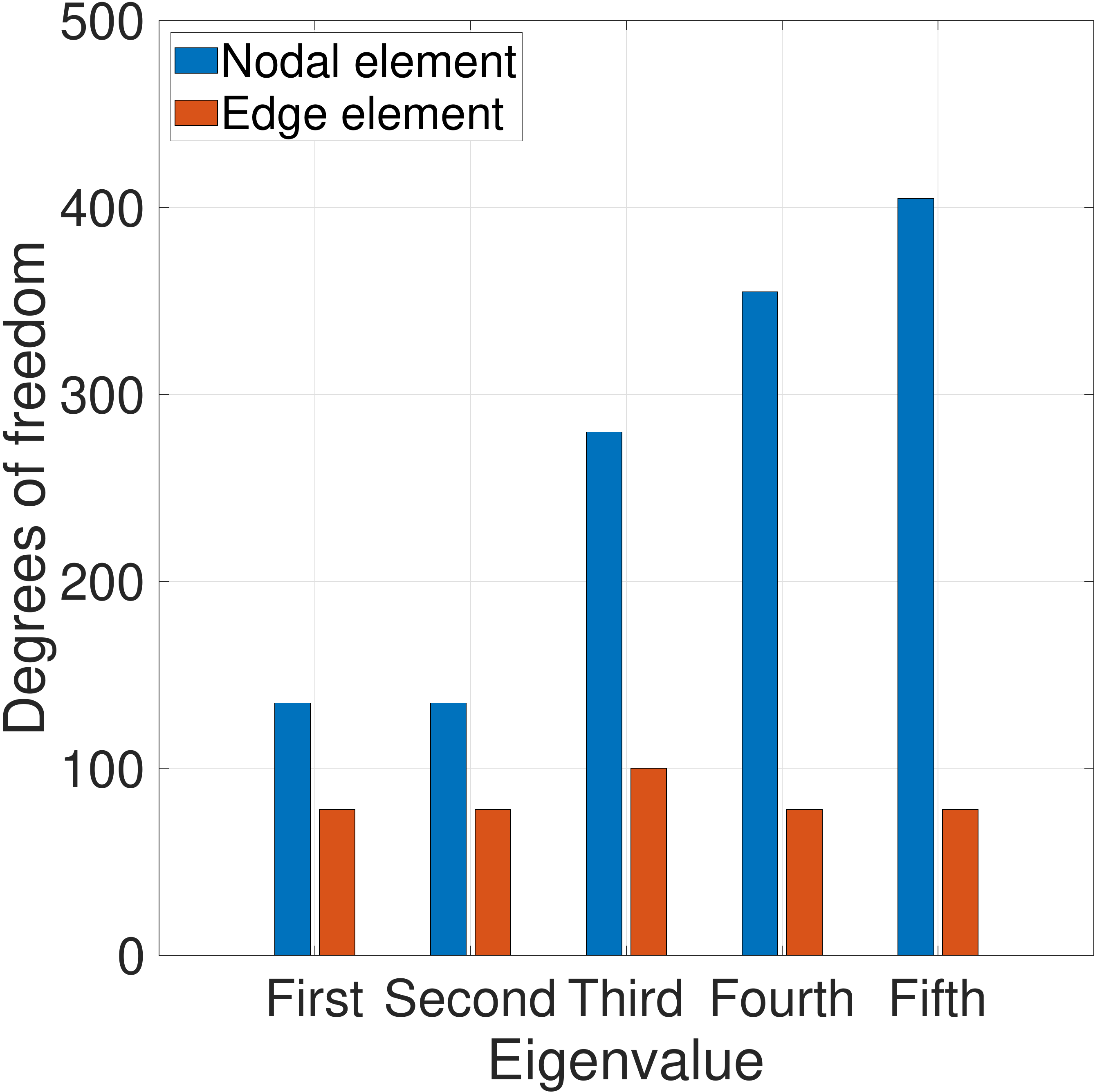}
\caption{Minimum number of FDOF required for less than 0.2$\%$ error (T6 vs ET8).}
\label{Sc}
\end{subfigure}
\caption{Comparative study of numerical performance of nodal elements with edge elements in predicting eigenvalues for square domain.}
\label{bar_square}
\end{figure}


The square of eigenvalues for each type of element is given in Table~\ref{tablesquare2} along with analytical results reported in~\cite{Boffi2001}. All of these elements have produced correct eigenvalues along with the correct multiplicities. For all the elements the first non zero eigenvalue has occurred after a certain number of zero eigen values at the machine precision level which signifies the approximation of null space. We have noticed a significant difference between nodal and edge elements in the mesh convergence analysis. With edge elements, we have found a far better coarse mesh accuracy. This fact is presented in bar diagrams in Fig.~\ref{bar_square}. We have compared Q9 and EQ12 in Fig.~\ref{Sa} on a scale of less than 1$\%$ error. Q4 and EQ4 elements are compared for less than 6$\%$ error in Fig.~\ref{Sb}. A scale of less than 0.2$\%$ error is chosen in Fig.~\ref{Sc} to compare T6 and ET8 elements. From Fig.~\ref{bar_square} we can see that we can attain the required level of accuracy for all five eigen values with a coarser edge element mesh than the required nodal element mesh for both quadrilateral and triangular elements.


\subsubsection{Circular shape domain}

\begin{table*}[pos=h!]
\caption{Analysis data of different nodal and edge elements for the circular shape domain problem.} \label{tabcircleshape1}
\begin{center}
\begin{tabular}{llllll}\hline
\multicolumn{3}{l}{{Nodal element}} &\multicolumn{3}{l}{{Edge element}} \\ \cmidrule(rr){1-3}\cmidrule(){4-6} 
\multicolumn{1}{l}{{Element}} & \multicolumn{1}{l}{{No. of}} & \multicolumn{1}{l}{{No. of}} & \multicolumn{1}{l}{{Element}} & \multicolumn{1}{l}{{No. of}} & \multicolumn{1}{l}{{No. of}} \\
\multicolumn{1}{l}{{Type}} & \multicolumn{1}{l}{{elements}} & \multicolumn{1}{l}{{FDOF/Equations}} & \multicolumn{1}{l}{{Type}} & \multicolumn{1}{l}{{elements}} & \multicolumn{1}{l}{{FDOF/Equations}}\\\cmidrule(rr){1-3} \cmidrule{4-6} 
Q9 $\&$ T6  & 256 & 2979 & EQ4 $\&$ ET3 & 1600 & 3160 \\
            &     &      & EQ12 $\&$ ET8 & 600 & 4720 \\  \hline
\end{tabular}
\end{center}
\end{table*}

\begin{figure}[pos=h!]
\centering
\includegraphics[width=0.7\columnwidth]{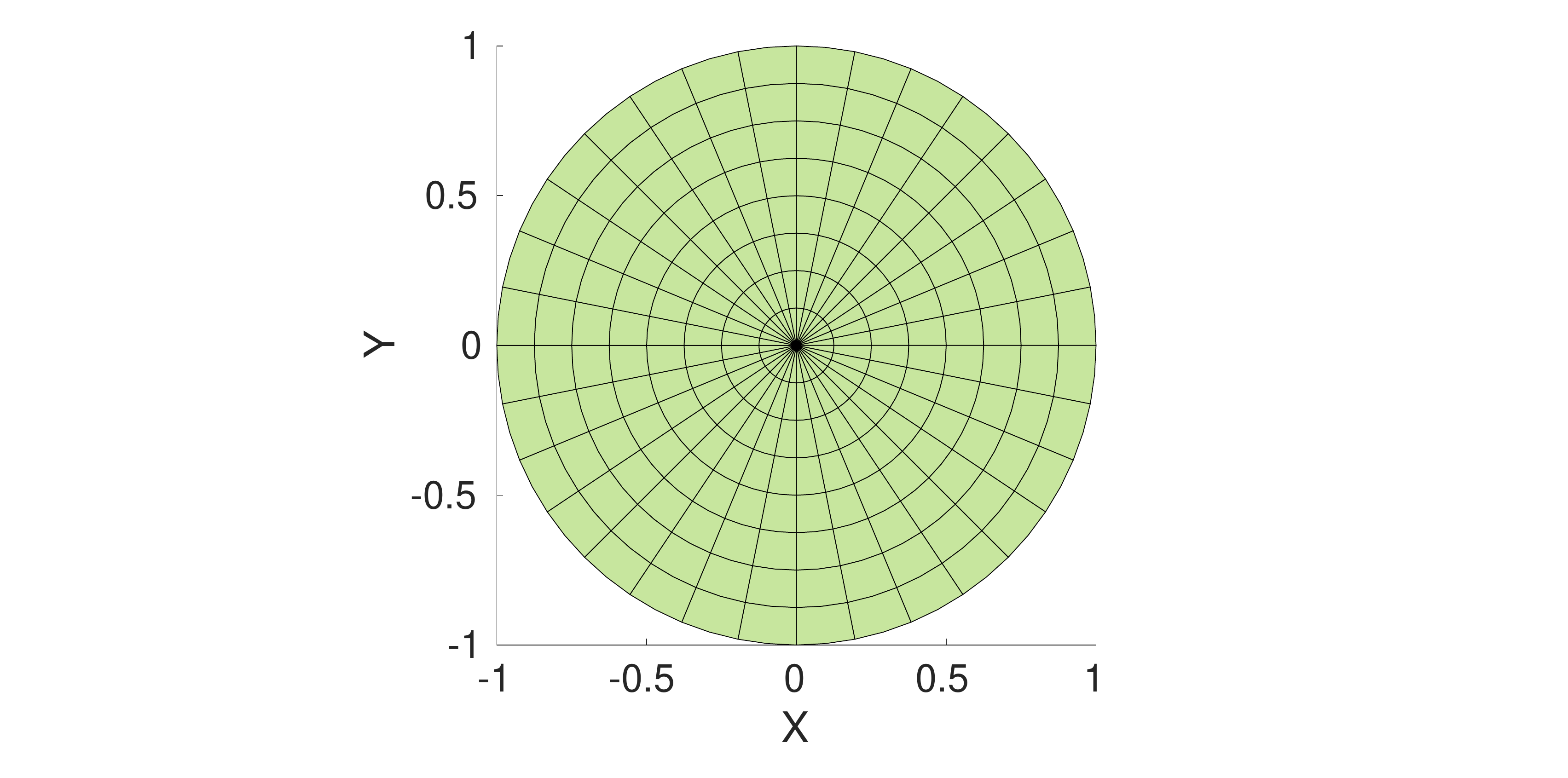}
\caption{Discretized circular domain}
\label{meshed_circle}
\end{figure}

\begin{table*}[pos=h!]
\begin{center}
\caption{$k_0^2$ on the circular domain for different elements (bracketed values show the multiplicity).} \label{tablecircle}
\begin{tabular}{llll}\hline
\multicolumn{1}{l}{Analytical} & \multicolumn{1}{l}{Nodal element} &\multicolumn{2}{l}{Edge element}  \\ \hline
Benchmark     &     Q9 $\&$ T6      &     EQ4 $\&$ ET3                       &      EQ12 $\&$ ET8          \\ \hline
3.391122 (2)  &  3.388867      &       3.410866                    &    3.380563            \\  
9.329970 (2)  &  9.318499      &       9.443712                    &    9.329356            \\
14.680392 (1) &  14.668347     &      14.768032                    &   14.756723            \\ 
17.652602 (2) &  17.609686     &      18.043890                    &   17.662206            \\
28.275806 (2) &  28.161350     &      28.586164                    &   28.262108            \\
28.419561 (2) &  28.376474     &      29.300835                    &   28.343555            \\
41.158640 (2) &  40.896578     &      43.394573                    &   41.404578            \\
44.970436 (2) &  44.858738     &      45.351002                    &   44.973322            \\
49.224256 (1) &  49.098151     &      49.598925                    &   49.656871            \\
56.272502 (2) &  55.747995     &      60.581974                    &   56.956837            \\
64.240225 (2) &  64.023638     &      65.163560                    &   64.271251            \\\hline
\multicolumn{3}{l}{Number of computed zeros} &  \multicolumn{1}{l}{}  \\  \midrule
-         &       992      &     401          &   450                 \\ \bottomrule
\end{tabular}
\end{center}
\end{table*}

\begin{figure}[pos=h!]
\centering
\includegraphics[width=0.7\columnwidth]{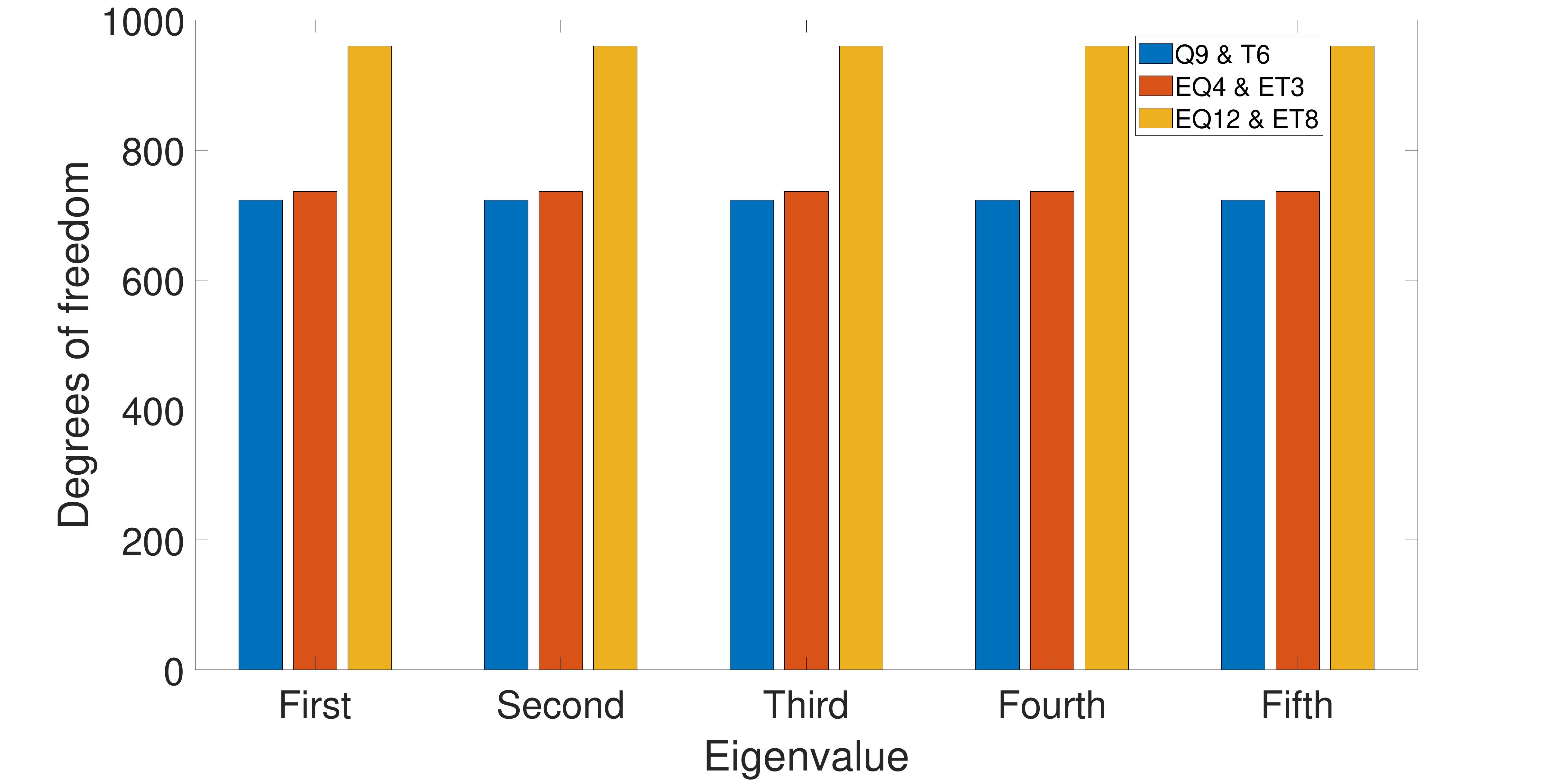}
\caption{Minimum number of FDOF required for less than 7$\%$ error in predicting eigen values for the circular domain}
\label{bar_circle}
\end{figure}

Here, the domain is a circle of radius unity and Fig.~\ref{meshed_circle} shows the discretized domain. The circular domain has a perfectly conducting boundary. No. of elements for different meshes of various element types are given in Table~\ref{tabcircleshape1}. The results are found with these meshes and they are closely matching with the benchmark values. For all the elements squared of the obtained eigenvalues are listed in Table~\ref{tablecircle}. The required degrees of freedom required to get less than 7$\%$ error for various elements has been presented in Fig.~\ref{bar_circle}. 

\subsubsection{L-shaped domain}

\begin{figure}[pos=h!]
\centering
\includegraphics[width=0.8\columnwidth]{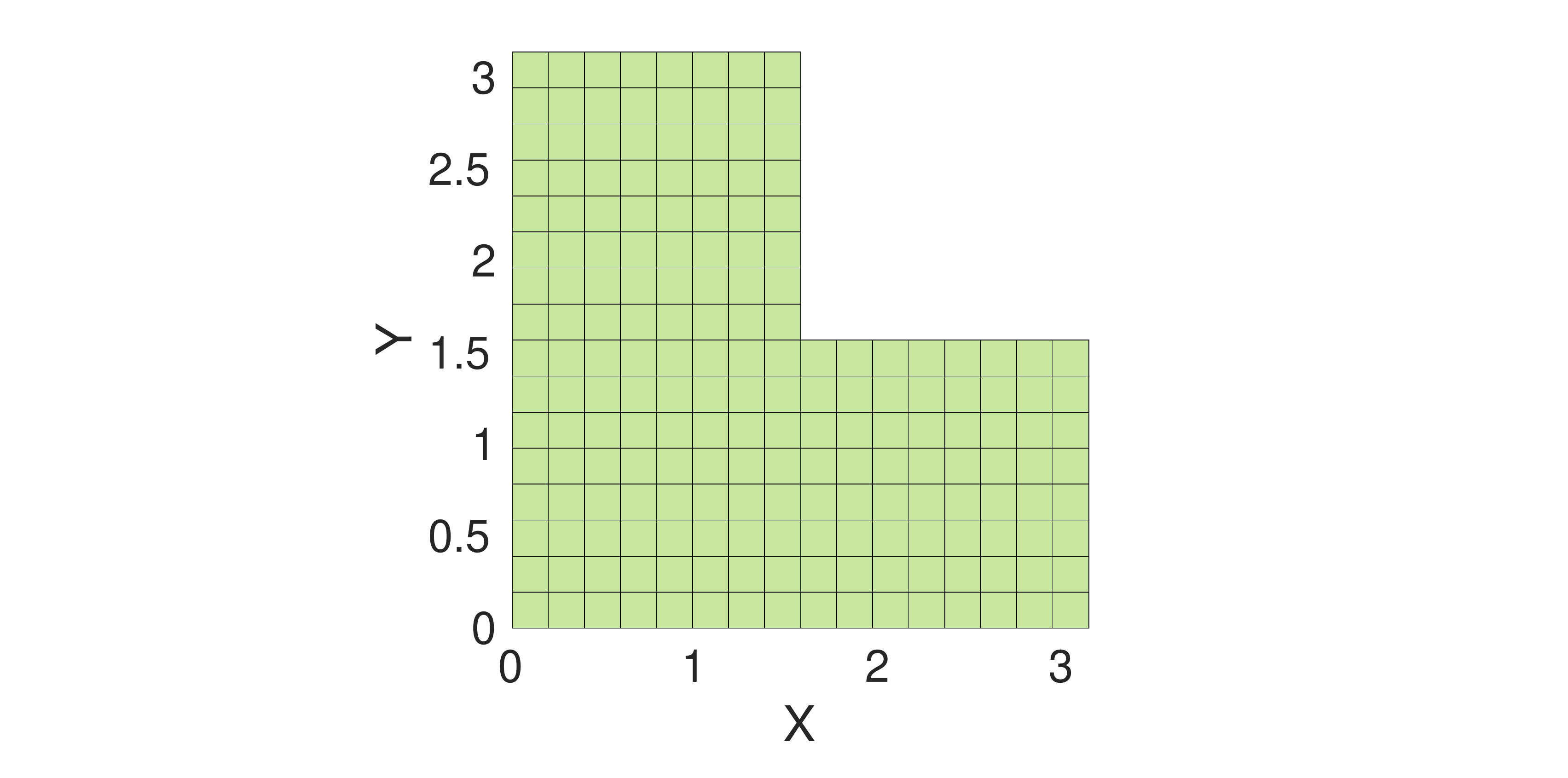}
\caption{Discretized L shape domain}
\label{meshed_Lshape}
\end{figure}
\begin{table*}[pos=h!]
\caption{Analysis data of different nodal and edge elements for the L shape domain problem.} \label{tabLshape1}
\begin{center}
\begin{tabular}{llllll}\hline
\multicolumn{3}{l}{{Nodal element}} &\multicolumn{3}{l}{{Edge element}} \\ \cmidrule(rr){1-3}\cmidrule(){4-6} 
\multicolumn{1}{l}{{Element}} & \multicolumn{1}{l}{{No. of}} & \multicolumn{1}{l}{{No. of}} & \multicolumn{1}{l}{{Element}} & \multicolumn{1}{l}{{No. of}} & \multicolumn{1}{l}{{No. of}} \\
\multicolumn{1}{l}{{Type}} & \multicolumn{1}{l}{{elements}} & \multicolumn{1}{l}{{FDOF/Equations}} & \multicolumn{1}{l}{{Type}} & \multicolumn{1}{l}{{elements}} & \multicolumn{1}{l}{{FDOF/Equations}}\\\cmidrule(rr){1-3} \cmidrule{4-6} 
Q4  & 768 & 2481 & EQ4  & 432 & 805 \\
T6  & 384 & 2481 & ET8  & 384 & 1856 \\
Q9  & 192 & 2481 & EQ12 & 192 & 1457 \\  \hline
\end{tabular}
\end{center}
\end{table*}
The L-shaped domain is obtained by deleting one quadrant from the square domain of side $\pi$ which has been considered in the previous example~(Section \ref{squaresection}). The discretized L shape domain with Q9 elements is shown in Fig.~\ref{meshed_Lshape}. The L-shaped domain is discretized with different mesh sizes for different elements, which is given in Table~\ref{tabLshape1}. 

\begin{figure}[pos=h!]
\centering
\begin{subfigure}{0.49\textwidth}
\centering
\includegraphics[width=0.8\textwidth]{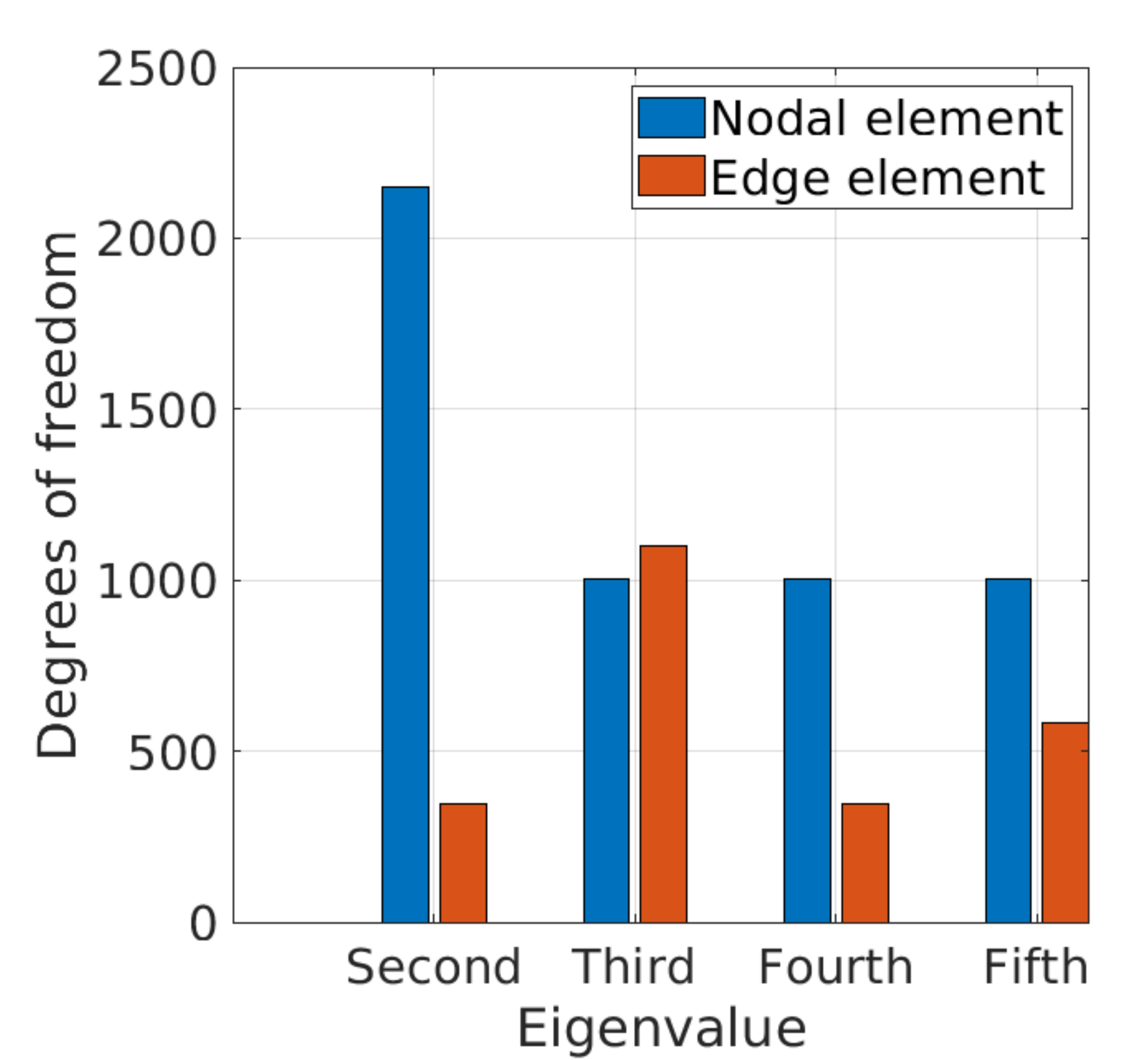}        
\caption{Minimum number of FDOF required for less than 6$\%$ error (Q9 vs EQ12).}
\label{Lshapea}
\end{subfigure}%
\hfill
\begin{subfigure}{0.49\textwidth}
\centering
\includegraphics[width=0.8\textwidth]{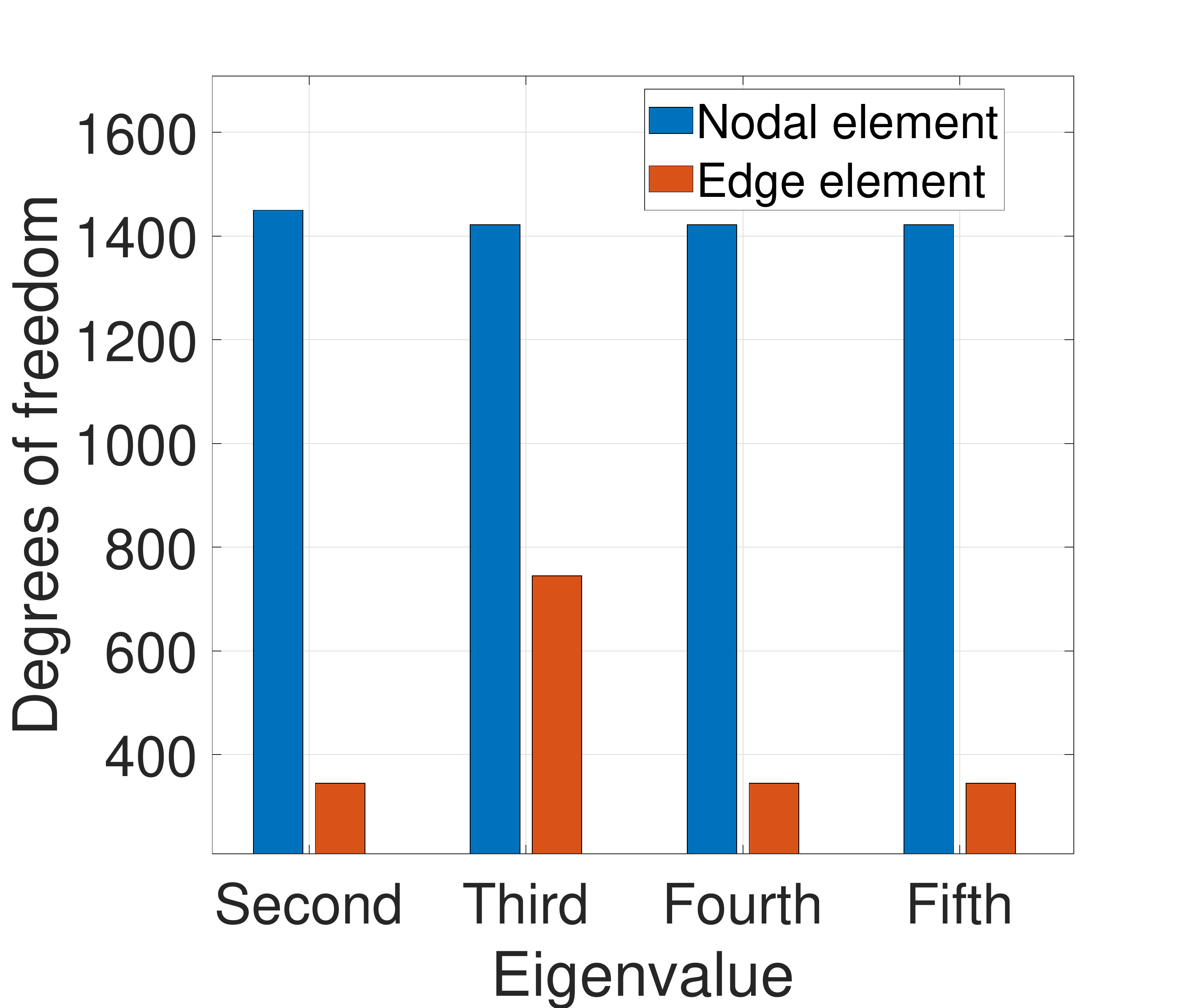}
\caption{Minimum number of FDOF required for less than 6$\%$ error (Q4 vs EQ4).}
\label{Lshapeb}
\end{subfigure}
\begin{subfigure}{0.49\textwidth}
\centering
\includegraphics[width=0.8\textwidth]{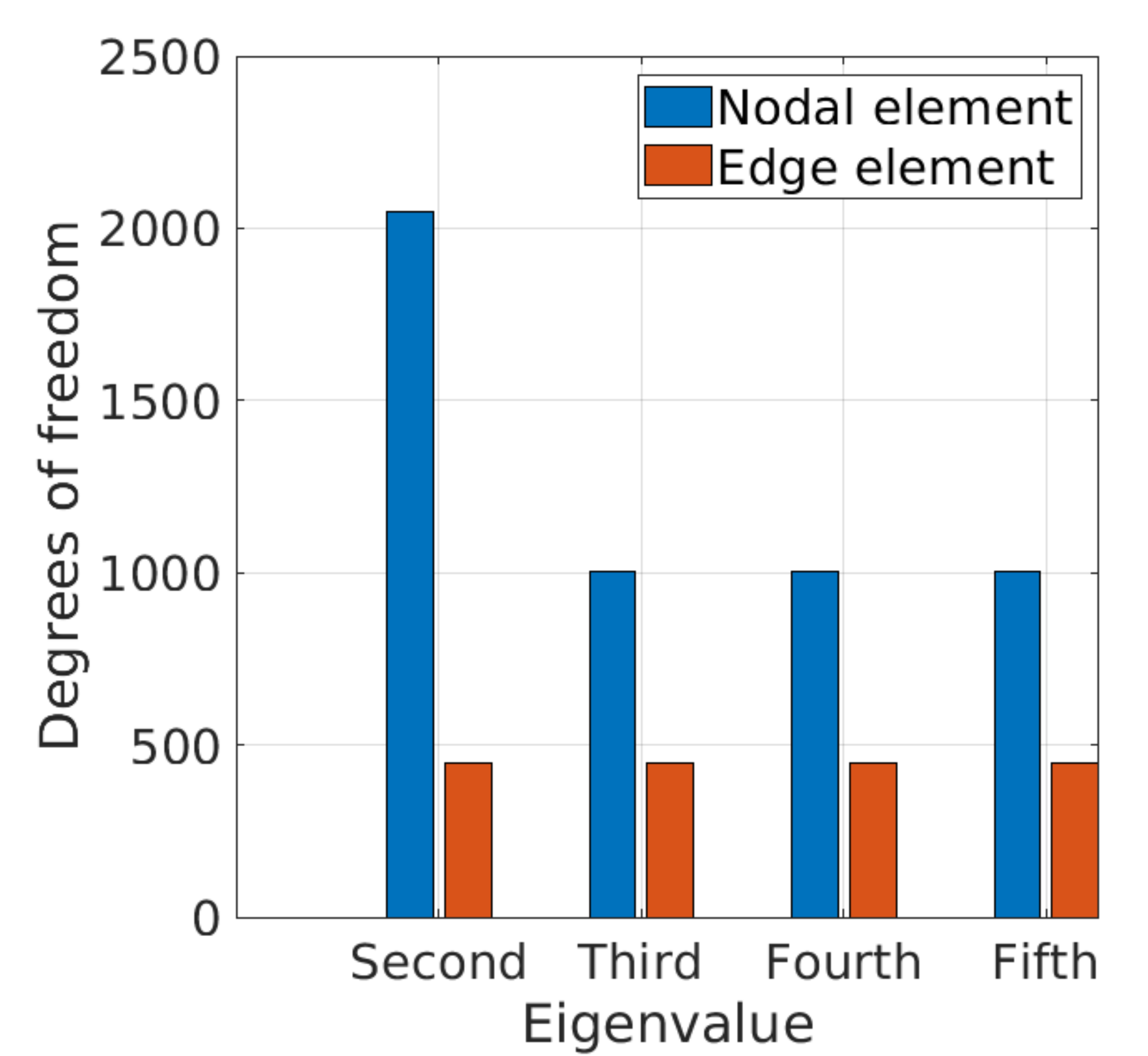}
\caption{Minimum number of FDOF required for less than 5$\%$ error (T6 vs ET8).}
\label{Lshapec}
\end{subfigure}
\caption{Comparative study of numerical performance of nodal elements with edge elements in predicting eigenvalues for L shape domain.}
\label{bar_lshape}
\end{figure}

As nodal elements can not capture the singular eigen value (0.591790), we have compared the accuracy of squared of the obtained eigenvalues in Fig.~\ref{bar_lshape} from the second eigen frequency. All the nodal elements are not able to capture the singular eigen value, as well as all of them, generate one spurious eigen value, whereas all the edge elements are able to predict the singular eigen value properly and they do not generate the spurious eigen value. We have found better coarse mesh accuracy with edge elements which is depicted in bar diagrams in Fig.~\ref{bar_lshape}. We have compared Q9 and EQ12 elements in Fig.~\ref{Lshapea} on a scale of less than 6$\%$ error. Q4 and EQ4 elements are compared on the scale of 6$\%$ error in Fig.~\ref{Lshapeb} whereas T6 and ET8 elements are compared for less than 4$\%$ error in Fig.~\ref{Lshapec}.
\begin{table*}[pos=h!]
\begin{center}
\caption{$k_0^2$ on the L-shaped domain for different elements.} \label{tabLshape2}
\begin{tabular}{lllllll}\hline
\multicolumn{1}{l}{} & \multicolumn{3}{l}{Nodal element}& \multicolumn{3}{l}{Edge element}  \\ \hline
  Benchmark  &       Q4     &     Q9      &     T6      &    EQ4     &   EQ12     &    ET8      \\ \hline
  0.591790   &       -      &     -       &     -       &  0.596170   & 0.597191   &  0.596538  \\
  1.432320   &    1.479654  &  1.507641   &  1.481174   &  1.434491   & 1.432148   &  1.432253  \\
     -       &    1.620830  &  2.277630   &  1.623752   &     -       &    -       &     -      \\
  4.005540   &    4.012869  &  3.997995   &  3.998716   &  3.781632   & 3.815079   &  3.999995  \\
  4.005540   &    4.012869  &  3.998150   &  3.998852   &  4.022899   & 4.000132   &  4.000019  \\
  4.613200   &    4.649897  &  4.645948   &  4.636011   &  4.418905   & 4.427041   &  4.616022  \\
  5.067330   &    5.577478  &  5.960861   &  5.561194   &  4.946604   & 4.939201   &  5.090472  \\
  7.955130   &    8.025738  &  7.990559   &  7.995340   &  7.821276   & 7.806319   &  8.000462  \\
  8.647370   &    9.404155  &  9.563233   &  9.336718   &  8.731623   & 8.649600   &  8.671775  \\
  9.481660   &    9.597178  &  9.838923   &  9.528452   &  9.575278   & 9.457106   &  9.460954  \\
 11.426100   &   12.380478  & 13.563361   &  12.304280  &  11.421820  & 11.349739  &  11.534996  \\
 14.448600   &   14.808709  & 14.771906   &  14.700390  &  14.611556  & 14.451696  &  14.542718  \\
 16.086200   &   16.206662  & 15.974127   &  15.984640  &  16.368790  & 16.008200  &  16.000365  \\ \hline
\multicolumn{4}{l}{Number of computed zeros}  &  \multicolumn{3}{l}{}       \\  \midrule
     -       &       828    &     827     &    827      &    64      &    39      &    22  \\ \bottomrule

\end{tabular}
\end{center}
\end{table*}

\vspace{5cm}
\subsubsection{Cracked circular domain}

\begin{figure}[pos=h!]
\centering
\includegraphics[width=0.3\textwidth]{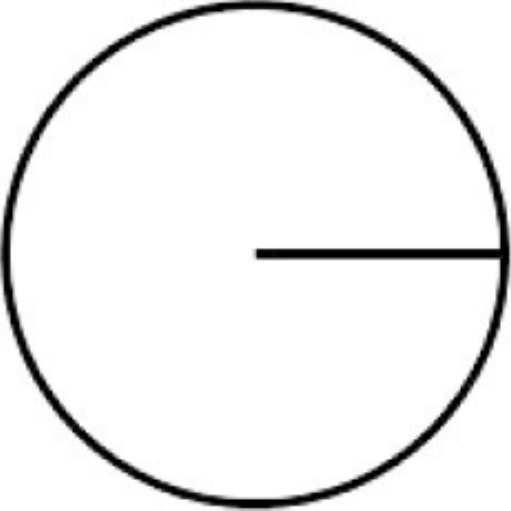}
\caption{Cracked Circular domain}
\label{crackedcircleshape}
\end{figure}
\begin{table*}[h!]
\begin{center}
\caption{Analysis data of different nodal and edge elements for the cracked circular domain problem.} \label{tabcrackcircleshape1}
\begin{tabular}{llllll}\hline
\multicolumn{3}{l}{{Nodal element}} &\multicolumn{3}{l}{{Edge element}} \\ \cmidrule(rr){1-3}\cmidrule(){4-6} 
\multicolumn{1}{l}{{Element}} & \multicolumn{1}{l}{{No. of}} & \multicolumn{1}{l}{{No. of}} & \multicolumn{1}{l}{{Element}} & \multicolumn{1}{l}{{No. of}} & \multicolumn{1}{l}{{No. of}} \\
\multicolumn{1}{l}{{Type}} & \multicolumn{1}{l}{{elements}} & \multicolumn{1}{l}{{FDOF/Equations}} & \multicolumn{1}{l}{{Type}} & \multicolumn{1}{l}{{elements}} & \multicolumn{1}{l}{{FDOF/Equations}}\\\cmidrule(rr){1-3} \cmidrule{4-6} 
Q9 $\&$ T6  & 256 & 3030 & EQ4 $\&$ ET3 & 384 & 723 \\
            &     &      & EQ12 $\&$ ET8 & 400 & 3161 \\  \hline
\end{tabular}
\end{center}
\end{table*}

\begin{table*}[h!]
\begin{center}
\caption{$k_0^2$ for the circular domain with crack for different elements.} \label{tablecrackedcircle}
\begin{tabular}{llll}\hline
\multicolumn{1}{l}{Analytical} & \multicolumn{1}{l}{Nodal element} &\multicolumn{2}{l}{Edge element}  \\ \hline
Benchmark     &     Q9 $\&$ T6      &     EQ4 $\&$ ET3      &      EQ12 $\&$ ET8          \\ \hline
1.358390      &      -         &   1.328769       &    1.243789            \\  
3.391122      &   3.732977     &   3.429584       &    3.368877            \\
6.059858      &   6.059539     &   6.155568       &    6.051724            \\ 
–             &   7.907027     &   -              &     -                  \\
9.329970      &   9.315916     &   9.540422       &    9.329118            \\
13.195056     &  13.171354     &  13.594326       &   13.201251            \\
14.680392     &  14.664092     &  14.941521       &   14.763221            \\
17.652602     &  17.601145     &  18.337851       &   17.662210            \\
21.196816     &  22.594381     &  21.136960       &   20.034173            \\
22.681406     &  28.140650     &  23.799963       &   22.709148            \\
28.275806     &  28.589967     &  28.986891       &   28.067093            \\\hline
\multicolumn{3}{l}{Number of computed zeros} &  \multicolumn{1}{l}{}  \\  \midrule
-         &       1009      &     40          &   8                   \\ \bottomrule
\end{tabular}
\end{center}
\end{table*}

Here, the domain is a circle of radius unity but it is with a crack as shown in Fig.~\ref{crackedcircleshape}.

\begin{figure}[pos=h!]
\centering
\includegraphics[width=0.6\textwidth]{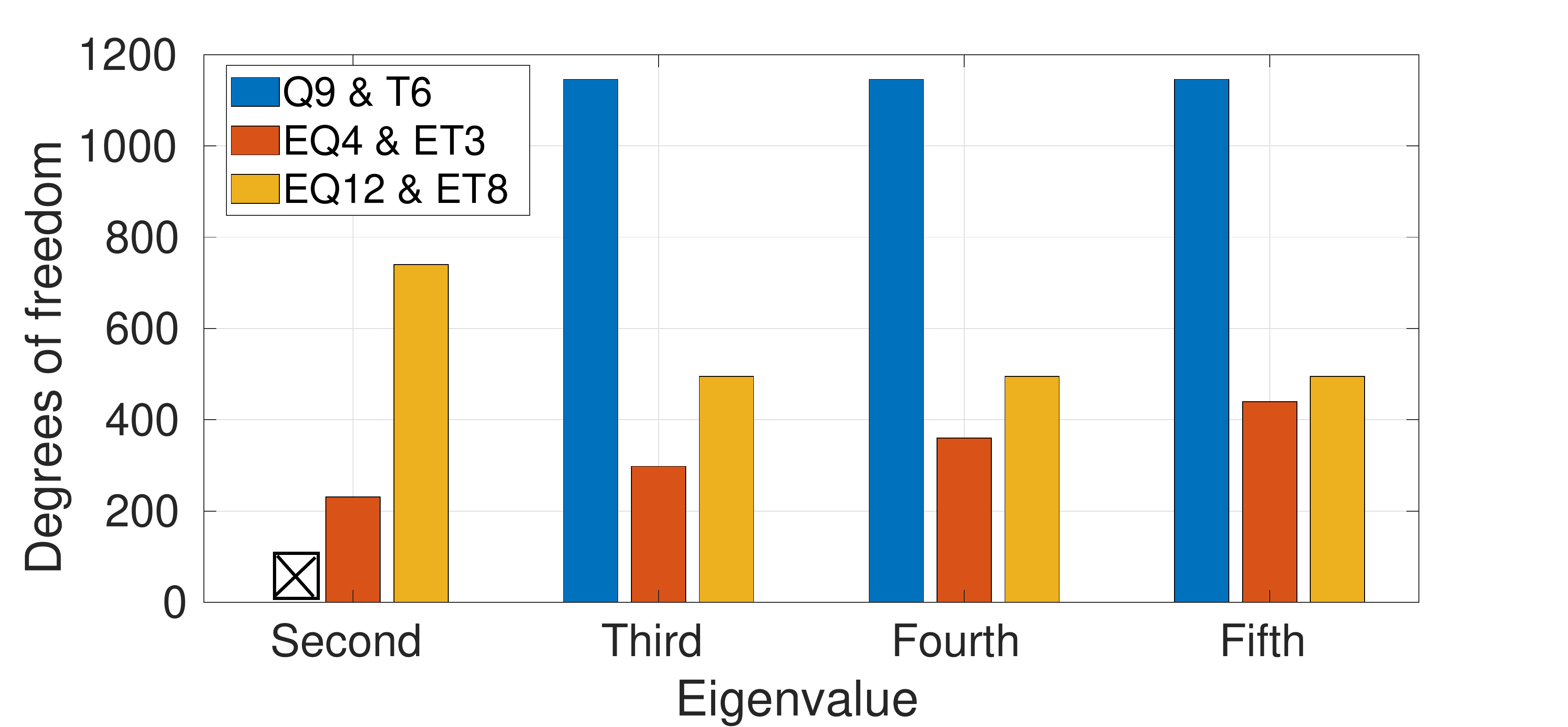}
\caption{Minimum number of FDOF required for less than 4$\%$ error in predicting eigen values for the cracked circular domain}
\label{bar_cracked_circle}
\end{figure}

The discretization of meshed domain with different types of elements is shown in Table~\ref{tabcrackcircleshape1}. The results are found with these meshes and they are closely matching with the benchmark values. As nodal elements can not capture the singular eigen value (1.358390), we have compared the accuracy of the obtained eigenvalues in Fig.~\ref{bar_cracked_circle} from the second eigen frequency. The degrees of freedom required to get an error less than 4$\%$ is shown in Fig.~\ref{bar_cracked_circle}. With Q9 $\&$ T6 elements, we have obtained error of 10.1$\%$ with 3030 FDOF for second eigen value.

\subsubsection{Curved L-shaped domain}   

\begin{figure}[h!]
\centering
\includegraphics[width=0.7\textwidth]{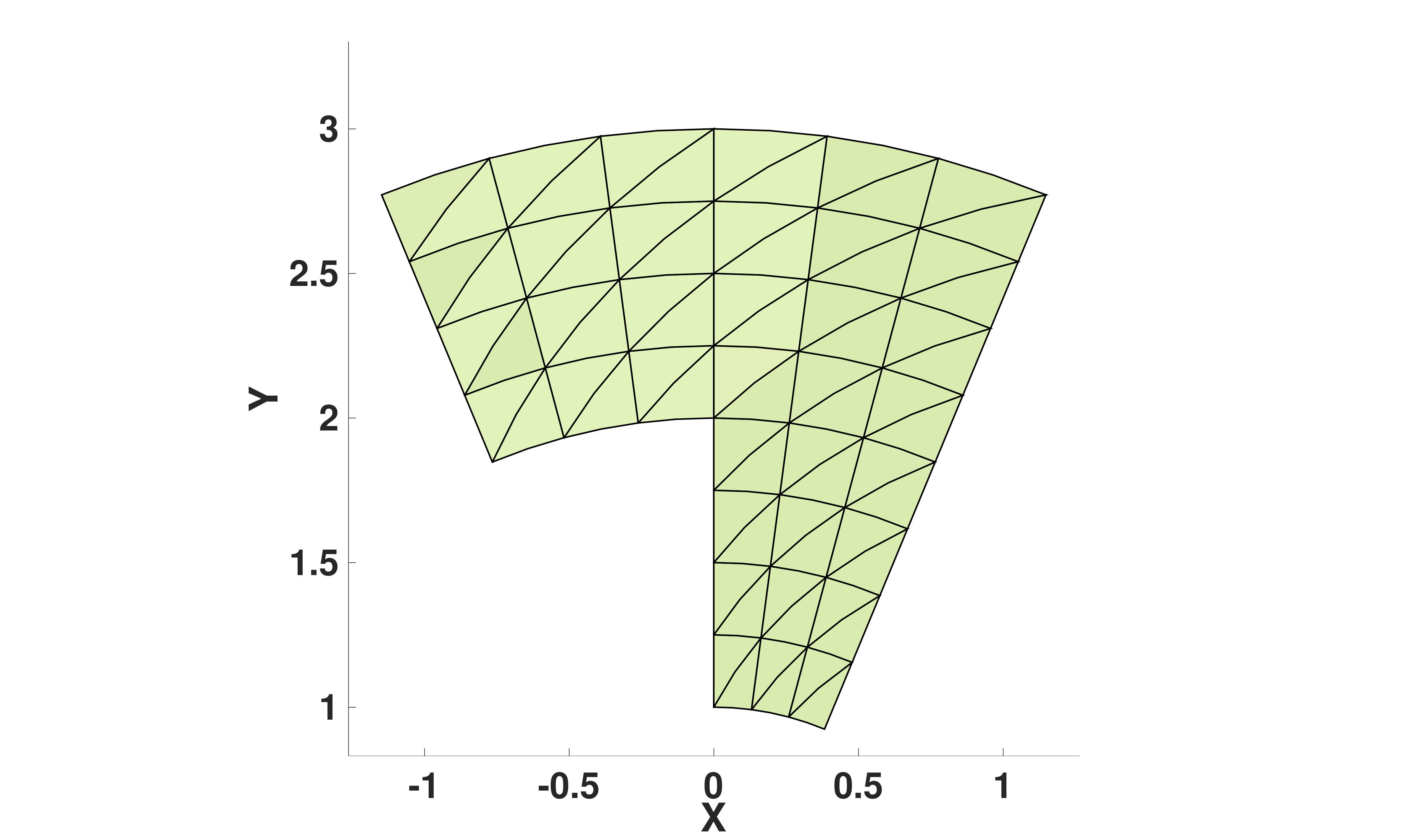}
\caption{Discretized curved L shape domain.}
\label{meshed_curvedl}
\end{figure}

\begin{table*}[pos=h!]
\caption{Analysis data of different nodal and edge elements for the curved L shape domain problem.} \label{tabshape1}
\begin{center}
\begin{tabular}{llllll}\hline
\multicolumn{3}{l}{{Nodal element}} &\multicolumn{3}{l}{{Edge element}} \\ \cmidrule(rr){1-3}\cmidrule(){4-6} 
\multicolumn{1}{l}{{Element}} & \multicolumn{1}{l}{{No. of}} & \multicolumn{1}{l}{{No. of}} & \multicolumn{1}{l}{{Element}} & \multicolumn{1}{l}{{No. of}} & \multicolumn{1}{l}{{No. of}} \\
\multicolumn{1}{l}{{Type}} & \multicolumn{1}{l}{{elements}} & \multicolumn{1}{l}{{FDOF/Equations}} & \multicolumn{1}{l}{{Type}} & \multicolumn{1}{l}{{elements}} & \multicolumn{1}{l}{{FDOF/Equations}}\\\cmidrule(rr){1-3} \cmidrule{4-6} 
Q4  & 270 & 909 & EQ4  & 200 & 480 \\
T6  & 96 & 657 & ET8  & 90 & 418 \\
Q9  & 48 & 657 & EQ12 & 84 & 790 \\  \hline
\end{tabular}
\end{center}
\end{table*}   

This example is taken from~\cite{benchmarkmonique}. Here, the domain has three straight and three circular sides of radii 1, 2, and 3, and Fig.~\ref{meshed_curvedl} shows the discretized domain.
 All the sides of the domain are perfectly conducting. This is one of the challenging problem as the domain is curved, non-convex along with sharp corner. Details of the meshed domain for different nodal and edge elements are given in Table~\ref{tabshape1}.

\begin{figure}[pos=h!]
\centering
\begin{subfigure}{0.49\textwidth}
\centering
\includegraphics[trim={10cm 0cm 10cm 0cm},width=0.8\textwidth]{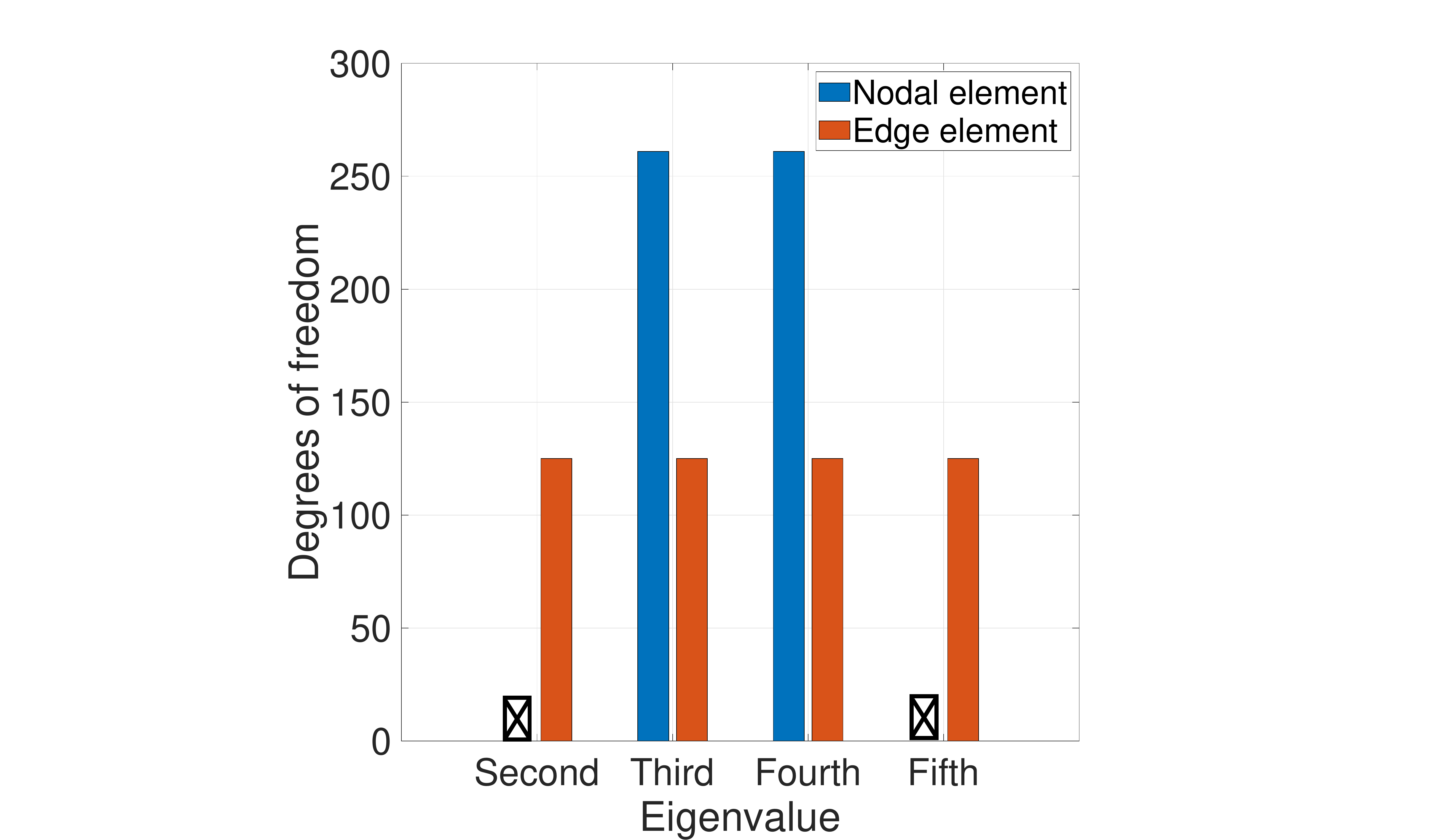}        
\caption{Minimum number of FDOF required for less than 2.5$\%$ error (Q9 vs EQ12).}
\label{La}
\end{subfigure}%
\hfill
\begin{subfigure}{0.49\textwidth}
\centering
\includegraphics[trim={10cm 0cm 9cm 0cm},width=0.8\textwidth]{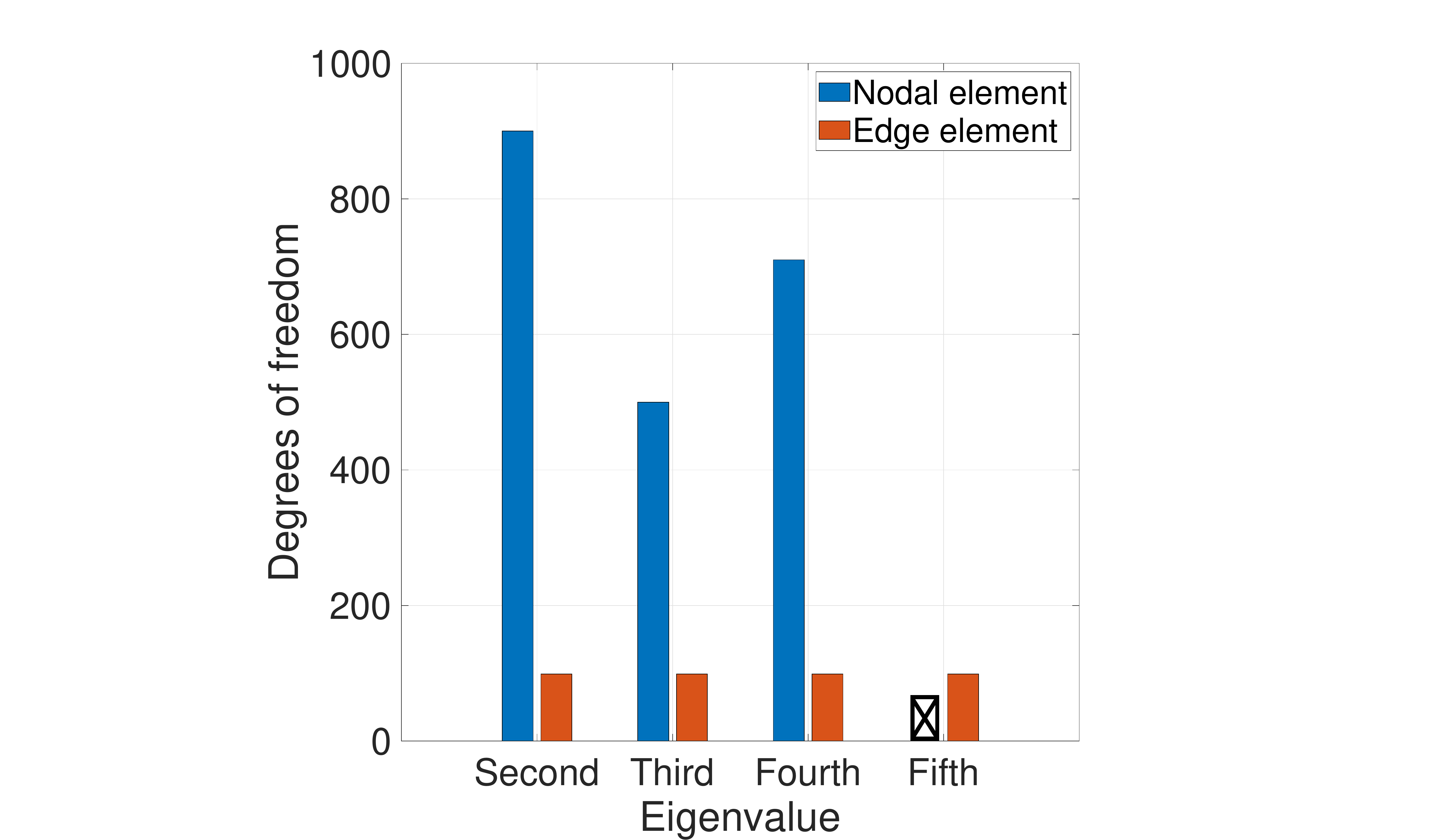}
\caption{Minimum number of FDOF required for less than 6.5$\%$ error (Q4 vs EQ4).}
\label{Lb}
\end{subfigure}
\begin{subfigure}{0.49\textwidth}
\centering
\includegraphics[trim={10cm 0cm 10cm 0cm},width=0.8\textwidth]{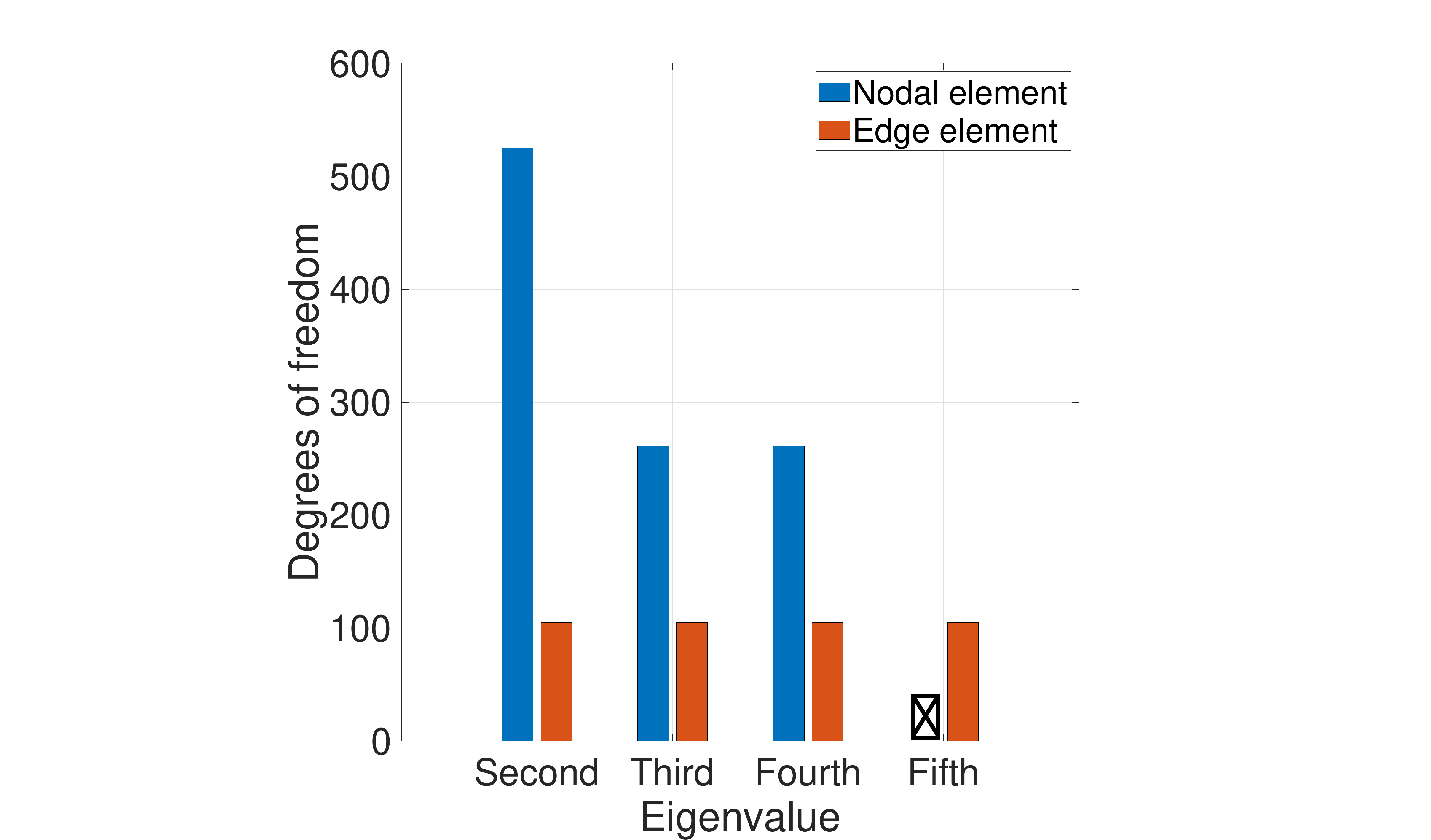}
\caption{Minimum number of FDOF required for less than 8$\%$ error (T6 vs ET8).}
\label{Lc}
\end{subfigure}
\caption{Comparative study of numerical performance of nodal elements with edge elements in predicting eigenvalues for curved L shape domain.}
\label{bar_curvedl}
\end{figure}

\begin{table*}[h!]
\begin{center}
\caption{$k_0^2$ on the curved L-shaped domain for different elements.} \label{tabshape2}
\begin{tabular}{lllllll}\hline
\multicolumn{1}{l}{} & \multicolumn{3}{l}{Nodal element} & \multicolumn{3}{l}{Edge element}  \\ \hline
  Benchmark  &       Q4    &     Q9      &     T6       &    EQ4     &   EQ12     &    ET8      \\ \hline
 1.818571   &        -     &     -       &      -      &  1.809583   & 1.814099   &  1.804651   \\
  3.490576   &    3.689097  &  3.854593   &  3.718116   &  3.505945   & 3.490535   &  3.489638  \\
     -       &    5.033585  &  6.889177   &  5.070643   &     -       &    -       &     -       \\
 10.065602   &   10.153471  & 10.047738   &  10.056494  &  10.198870  & 10.066721  &  10.071826  \\
 10.111886   &   10.284415  & 10.252722   &  10.194525  &  10.241690  & 10.116060  &  10.109109  \\
 12.435537   &   13.845497  & 15.141449   &  13.759140  &  12.554637  & 12.427678  &  12.420649  \\ \hline
\multicolumn{4}{l}{Number of computed zeros}  &  \multicolumn{3}{l}{}       \\  \midrule
     -       &       303     &     219     &    219      &    184      &    173      &    149 \\ \bottomrule

\end{tabular}
\end{center}
\end{table*}
 For all the elements, squared of the obtained eigenvalues are listed in Table~\ref{tabshape2} along with the benchmark values from~\cite{benchmarkmonique}. All the nodal elements are not able to capture the singular eigen value (1.818571) (due to the presence of a sharp corner) as well as all of them generate one spurious eigen value. Whereas all the edge elements are able to predict the singular eigen value properly and they do not generate the spurious eigen value. Therefore we have compared accuracy from the second eigenvalue in Fig.~\ref{bar_curvedl}.  We have found better coarse mesh accuracy with edge elements which is depicted in bar diagrams in Fig.~\ref{bar_curvedl}. We have compared Q9 and EQ12 elements in Fig.~\ref{La} on a scale of less than 2.5$\%$ error. A scale of 6.5$\%$ error is chosen in Fig.~\ref{Lb} to compare Q4 and EQ4. T6 and ET8 elements are compared for less than 8$\%$ error in Fig.~\ref{Lc}. With 657 FDOF of Q9 element, we have obtained 10.4$\%$ error for the second eigen value and 21.7$\%$ error for the fifth eigen value. An error of 11.3$\%$ is obtained for Q4 element with 909 FDOF for the fifth eigen value. With T6 element for the fifth eigen value we get an error of 10.6$\%$ with 657 FDOF.

\subsubsection{Inhomogeneous L shape domain} 

\begin{figure}[pos=h!]
\centering
\includegraphics[width=0.8\textwidth]{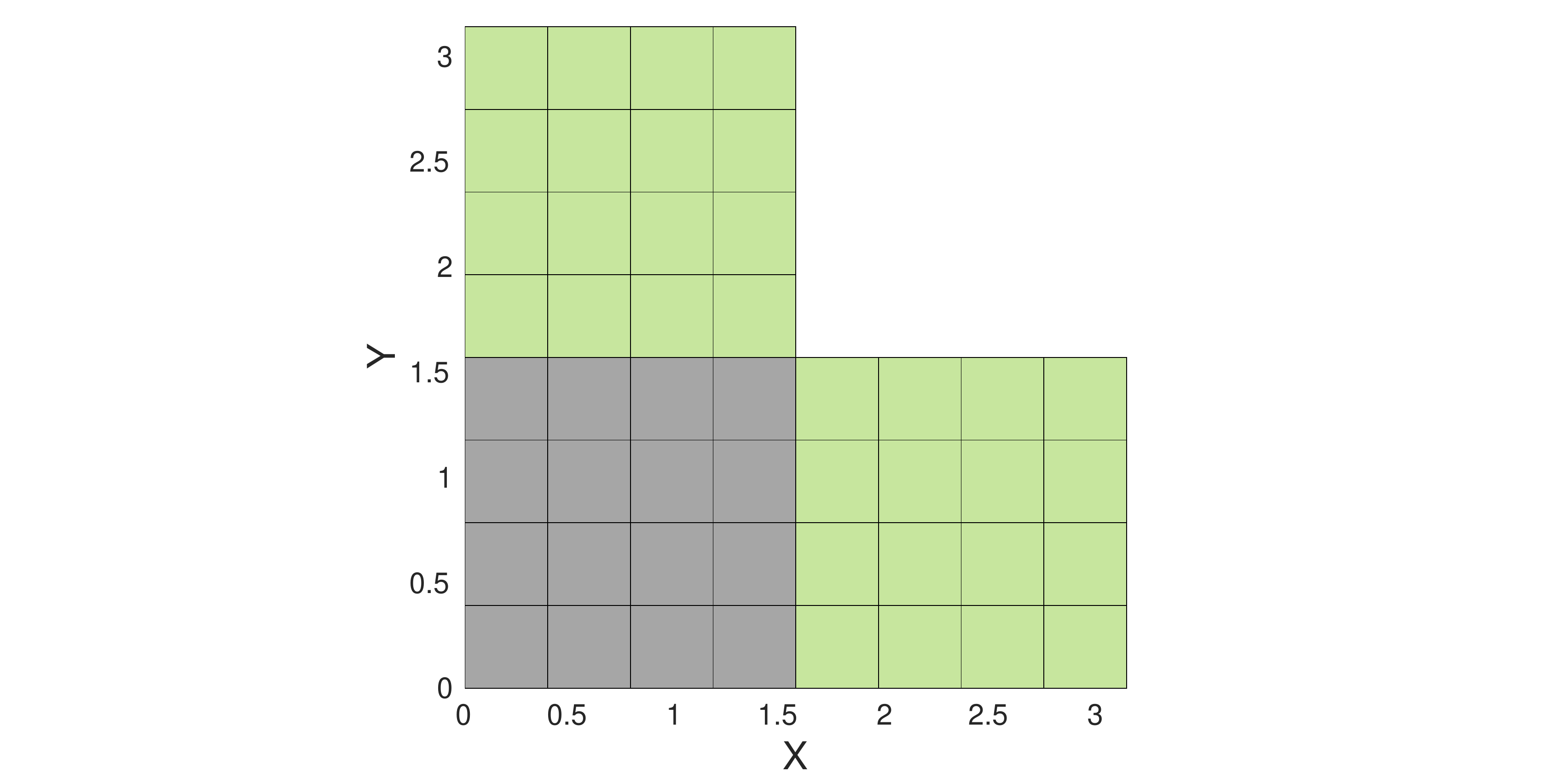}
\caption{Inhomogeneous L shape domain.}
\label{inhomogeneousL}
\end{figure}

\begin{table*}[pos=h!]
\begin{center}
\caption{Analysis data of different nodal and edge elements for the inhomogeneous L shape domain problem.} \label{tabLshape3}
\begin{tabular}{llllll}\hline
\multicolumn{3}{l}{{Nodal element}} &\multicolumn{3}{l}{{Edge element}} \\ \cmidrule(rr){1-3}\cmidrule(){4-6} 
\multicolumn{1}{l}{{Element}} & \multicolumn{1}{l}{{No. of}} & \multicolumn{1}{l}{{No. of}} & \multicolumn{1}{l}{{Element}} & \multicolumn{1}{l}{{No. of}} & \multicolumn{1}{l}{{No. of}} \\
\multicolumn{1}{l}{{Type}} & \multicolumn{1}{l}{{elements}} & \multicolumn{1}{l}{{FDOF/Equations}} & \multicolumn{1}{l}{{Type}} & \multicolumn{1}{l}{{elements}} & \multicolumn{1}{l}{{FDOF/Equations}}\\\cmidrule(rr){1-3} \cmidrule{4-6} 
Q4  & 768 & 2481 & EQ4  & 768 & 1457 \\
T6  & 384 & 2481 & ET8  & 384 & 1856 \\
Q9  & 192 & 2481 & EQ12 & 192 & 1457 \\  \hline
\end{tabular}
\end{center}
\end{table*}

\begin{table*}[pos=h!]
\begin{center}
\caption{$k_0^2$ on the inhomogeneous L shape domain for different elements.} \label{tabLshape4}
\begin{tabular}{lllllll}\hline
\multicolumn{1}{l}{} & \multicolumn{3}{l}{Nodal element}& \multicolumn{3}{l}{Edge element}  \\ \hline
  Benchmark  &       Q4     &     Q9      &     T6      &    EQ4     &   EQ12     &    ET8      \\ \hline
  0.175980   &       -      &     -       &     -       &  0.176192   & 0.176090   &  0.176085  \\
  0.398080   &    0.410869  &  0.418399   &  0.411088   &  0.397985   & 0.397393   &  0.397395  \\
     -       &    0.410993  &  0.526242   &  0.411234   &     -       &    -       &     -      \\
  0.964840   &    0.973806  &  0.973085   &  0.970759   &  0.872179   & 0.853464   &  0.966204  \\
  0.978740   &    0.998372  &  1.003343   &  0.994894   &  0.976496   & 0.973764   &  0.980699  \\
  1.524310   &    1.785034  &  1.783645   &  1.775821   &  1.529238   & 1.521521   &  1.522974  \\
  1.765930   &    1.816265  &  1.973428   &  1.807250   &  1.768005   & 1.758150   &  1.761401  \\
  2.274180   &    2.306054  &  2.292383   &  2.293151   &  2.252808   & 2.233477   &  2.293116  \\
  2.389530   &    2.563009  &  2.694188   &  2.544549   &  2.395908   & 2.381229   &  2.412659  \\
  3.394090   &    3.428865  &  3.380034   &  3.382655   &  3.424851   & 3.385157   &  3.384739  \\
  3.397400   &    3.433245  &  3.384257   &  3.386683   &  3.427972   & 3.388137   &  3.388719  \\
  3.646940   &    3.695923  &  3.664321   &  3.663456   &  3.664173   & 3.630283   &  3.641584  \\
  3.664270   &    3.906853  &  4.096015   &  3.873663   &  3.686256   & 3.655079   &  3.658646  \\ \hline
\multicolumn{4}{l}{Number of computed zeros}  &  \multicolumn{3}{l}{}       \\  \midrule
     -       &       827    &     827     &    827      &    19      &    56      &    15  \\ \bottomrule
\end{tabular}
\end{center}
\end{table*}
In this example, the efficacy of the edge elements for the inhomogeneous domain becomes evident. Fig~\ref{inhomogeneousL} represents the domain where relative permittivities are 1 and 5 for grey and the light green regions respectively. Discretization details of the domain for various elements are presented in Table~\ref{tabLshape3}. For all the elements squared of the obtained eigenvalues are listed in Table.~\ref{tabLshape4}. All the nodal elements are not able to capture the singular eigen value (0.175980) as well as all of them generate one spurious eigen value.

\begin{figure}[pos=h!]
\centering
\begin{subfigure}{0.49\textwidth}
\centering
\includegraphics[width=0.8\textwidth]{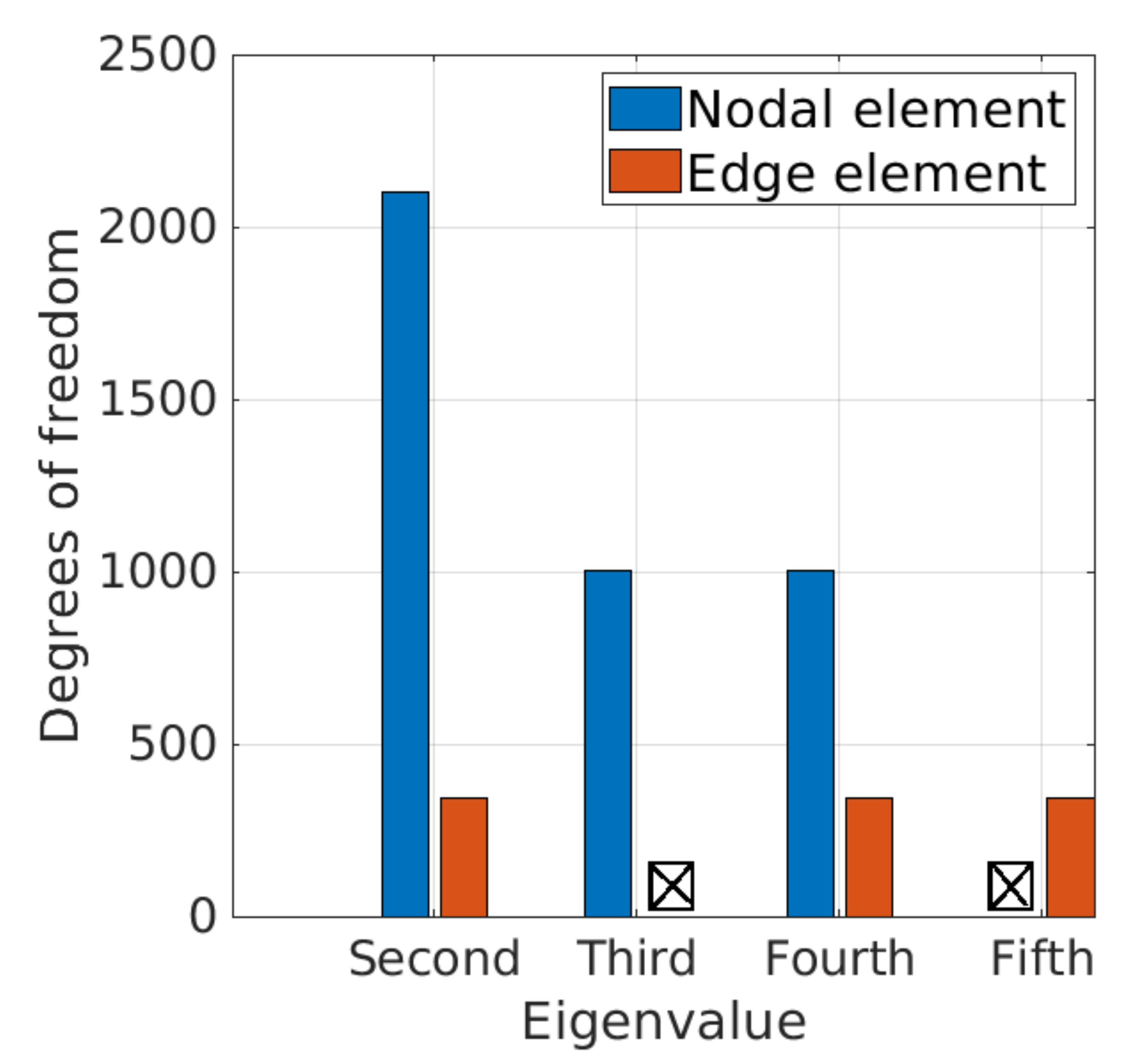}        
\caption{Minimum number of FDOF required for less than 6$\%$ error (Q9 vs EQ12).}
\label{Lshapea2}
\end{subfigure}%
\hfill
\begin{subfigure}{0.49\textwidth}
\centering
\includegraphics[width=0.8\textwidth]{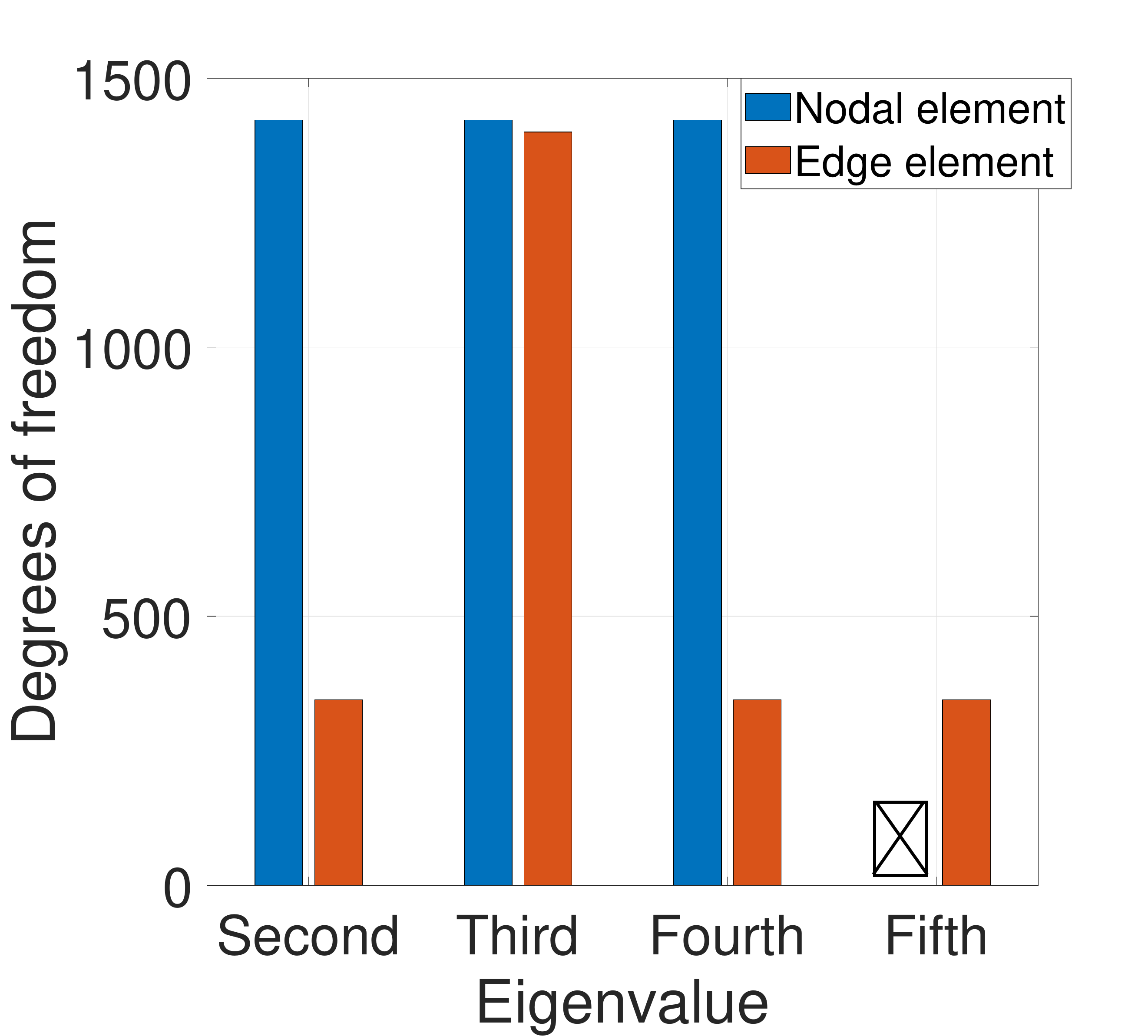}
\caption{Minimum number of FDOF required for less than 10$\%$ error (Q4 vs EQ4).}
\label{Lshapeb2}
\end{subfigure}
\begin{subfigure}{0.49\textwidth}
\centering
\includegraphics[width=0.8\textwidth]{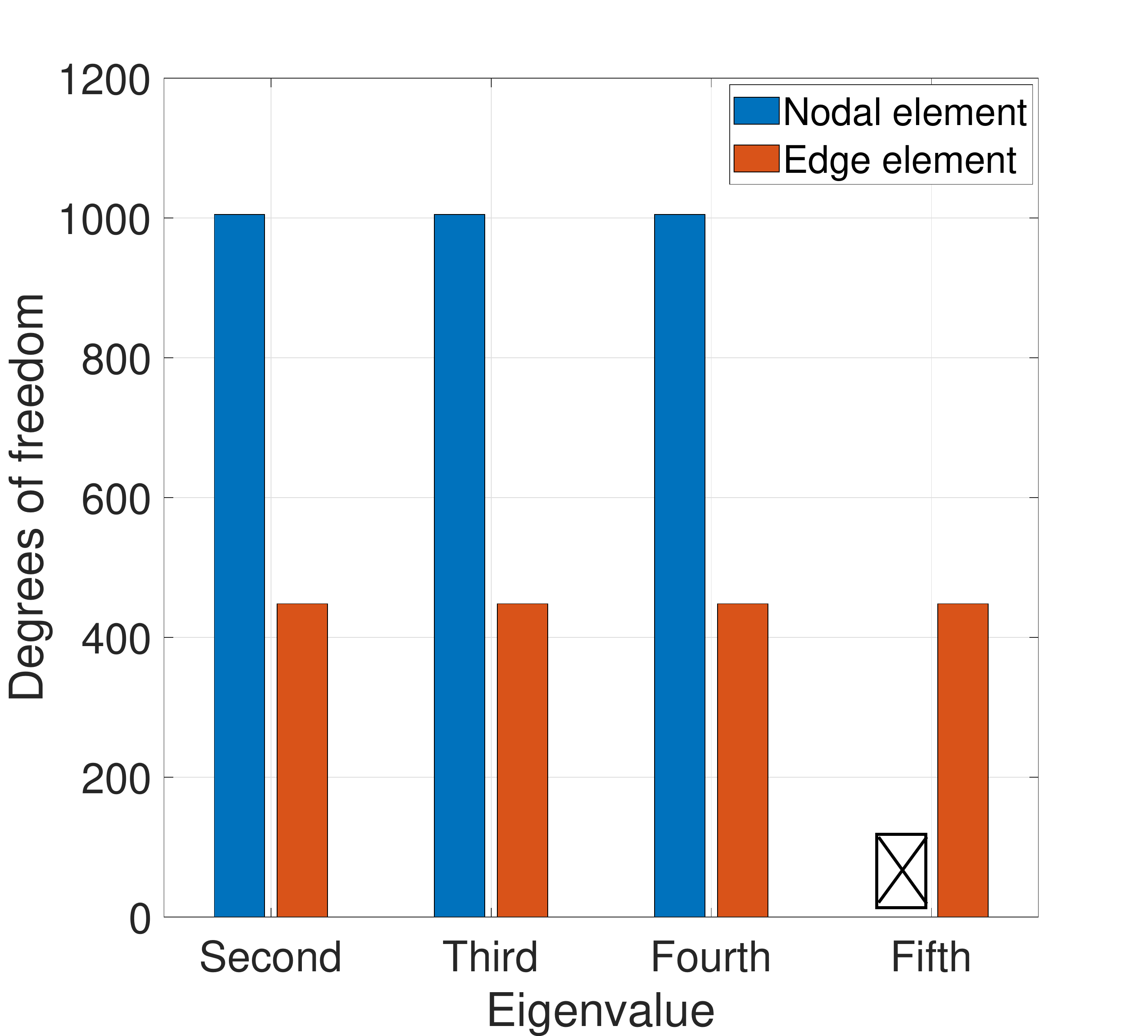}
\caption{Minimum number of FDOF required for less than 5$\%$ error (T6 vs ET8).}
\label{Lshapec2}
\end{subfigure}
\caption{Comparative study of numerical performance of nodal elements with edge elements in predicting eigenvalues for inhomogeneous L shape domain.}
\label{bar_lshape2}
\end{figure}

Whereas all the edge elements are able to predict the singular eigen value properly and they do not generate the spurious eigen value, thus we have compared the accuracy of elements from the second eigen value in Fig.~\ref{bar_lshape2}.We have found better coarse mesh accuracy with edge elements which is depicted in bar diagrams in Fig.~\ref{bar_lshape2}. We have compared Q9 and EQ12 elements in Fig.~\ref{Lshapea2} on a scale of less than 6$\%$ error. Q4 and EQ4 elements are compared on a scale of 10$\%$ error in Fig.~\ref{Lshapeb2} whereas T6 and ET8 elements are compared for less than 5$\%$ error in Fig.~\ref{Lshapec2}.

With Q9 elements for 2481 FDOF there is 17$\%$ error for the fifth eigen value. For EQ12 elements we obtained error of 11.8$\%$ with 1457 FDOF for the third eigen value while for Q4 element we get 17.1$\%$ error for 2481 FDOF with the fifth eigen value. Error of 16.5$\%$ is obtained with T6 element for 2481 FDOF for the fifth eigen value.

\vspace{4cm}
\subsection{Performance analysis for distorted mesh}


\begin{table}[pos=h!]
\begin{center}
\caption{$k_0^2$ on the curved L shape domain for normal and distorted meshes of lower order higher order of nodal and edge elements.} \label{tabcurlshape} 
\begin{tabular}{p{1.9cm}p{1.9cm}p{2.0cm}p{1.8cm}p{1.8cm}}   
 \multicolumn{5}{c}{} \\  
\multicolumn{3}{l}{Q4: 198 elements (681 FDOF) } &\multicolumn{2}{l}{EQ4: 198 elements (362 FDOF) } \\  \hline\hline
\multicolumn{1}{l}{Analytical} & \multicolumn{2}{l}{Nodal element (Q4)}  &\multicolumn{2}{l}{Edge element (EQ4)} \\\cmidrule(ll){1-3}\cmidrule(){4-5}
    -      & Normal    & Distorted &    Normal   & Distorted  \\\cmidrule(rr){1-3}\cmidrule(){4-5}
1.818571   &    -      &    -      &    1.810926 &  1.810926  \\
3.490576   & 3.725217  & 3.729079  &    3.508771 &  3.508770  \\
     -     & 5.084991  & 5.085667  &        -    &      -     \\        
10.065602  & 10.145899 & 10.148647 &   10.141604 &  10.141603 \\
10.111886  & 10.386313 & 10.407295 &   10.348577 &  10.348577 \\
12.435537  & 13.897028 & 13.902881 &   12.501734 &  12.501734 \\\cmidrule(rr){1-3} \cmidrule(){4-5} 
\multicolumn{4}{l}{Number of computed zeros}  &  \multicolumn{1}{l}{}       \\ \cmidrule(rr){1-3}\cmidrule{4-5}
	-       &    65       &   65    &   227     &   227            \\ \hline 
 \multicolumn{5}{c}{} \\ 
\multicolumn{3}{l}{Q9: 75 elements (1005 FDOF) } &\multicolumn{2}{l}{EQ12: 75 elements (560 FDOF)} \\  \hline\hline 
\multicolumn{1}{l}{Analytical} & \multicolumn{2}{l}{Nodal element (Q9)}  &\multicolumn{2}{l}{Edge element (EQ12)} \\ \cmidrule(ll){1-3}\cmidrule(){4-5}
-          &   Normal  & Distorted  &    Normal     & Distorted   \\\cmidrule(rr){1-3}\cmidrule{4-5}
1.818571   &    -      &      -     &   1.813849    &   1.816383  \\
3.490576   & 3.798913  &  3.799074  &   3.490505    &   3.495839  \\
     -     & 6.828452  &  6.843561  &         -     &       -     \\        
10.065602  & 10.049142 & 10.058195  &   10.067992   &   10.088585 \\
10.111886  & 10.242052 & 10.245292  &   10.113374   &   10.140025 \\
12.435537  & 15.071954 & 15.085305  &   12.429564   &   12.443714 \\ \cmidrule(rr){1-3}\cmidrule(){4-5}
\multicolumn{3}{l}{Number of computed zeros}  &  \multicolumn{2}{l}{}       \\ \cmidrule{1-3}\cmidrule(){4-5}
   -       &    335    &    335     &    261        &   261            \\ \hline
 \multicolumn{5}{c}{} \\ 
\multicolumn{3}{l}{T6: 72 elements (501 FDOF)} &\multicolumn{2}{l}{ET8: 72 elements (332 FDOF)} \\ \hline\hline 
\multicolumn{1}{l}{Analytical} & \multicolumn{2}{l}{Nodal element (Q9)} &\multicolumn{2}{l}{Edge element (EQ12)} \\\cmidrule(ll){1-3}\cmidrule(){4-5}
-          &   Normal  & Distorted &  Normal     & Distorted  \\\cmidrule(rr){1-3}\cmidrule{4-5}
1.818571   &    -      &      -    & 1.803444    &  1.799088  \\
3.490576   &3.753233   & 3.751064  & 3.489667    &  3.484439  \\
     -     &5.118928   & 5.109720  &      -      &         -  \\        
10.065602  &10.063793  & 10.063864 &10.070665    &  10.059423 \\
10.111886  &10.204527  & 10.212374 &10.110606    &  10.110596 \\
12.435537  &13.806202  & 13.824826 &12.414018    &  12.407989 \\ \cmidrule(rr){1-3}\cmidrule{4-5}
\multicolumn{3}{l}{Number of computed zeros}  &  \multicolumn{2}{l}{}       \\ \cmidrule(rr){1-3}\cmidrule{4-5}
   -       &    167    &    167    &   117        &   117            \\ \hline\hline
\end{tabular}
\end{center}
\end{table}

\begin{figure}[pos=h!]
\centering
\begin{subfigure}{0.3\columnwidth}
\centering
\includegraphics[trim={16cm 0cm 16cm 0cm},width=0.48\columnwidth]{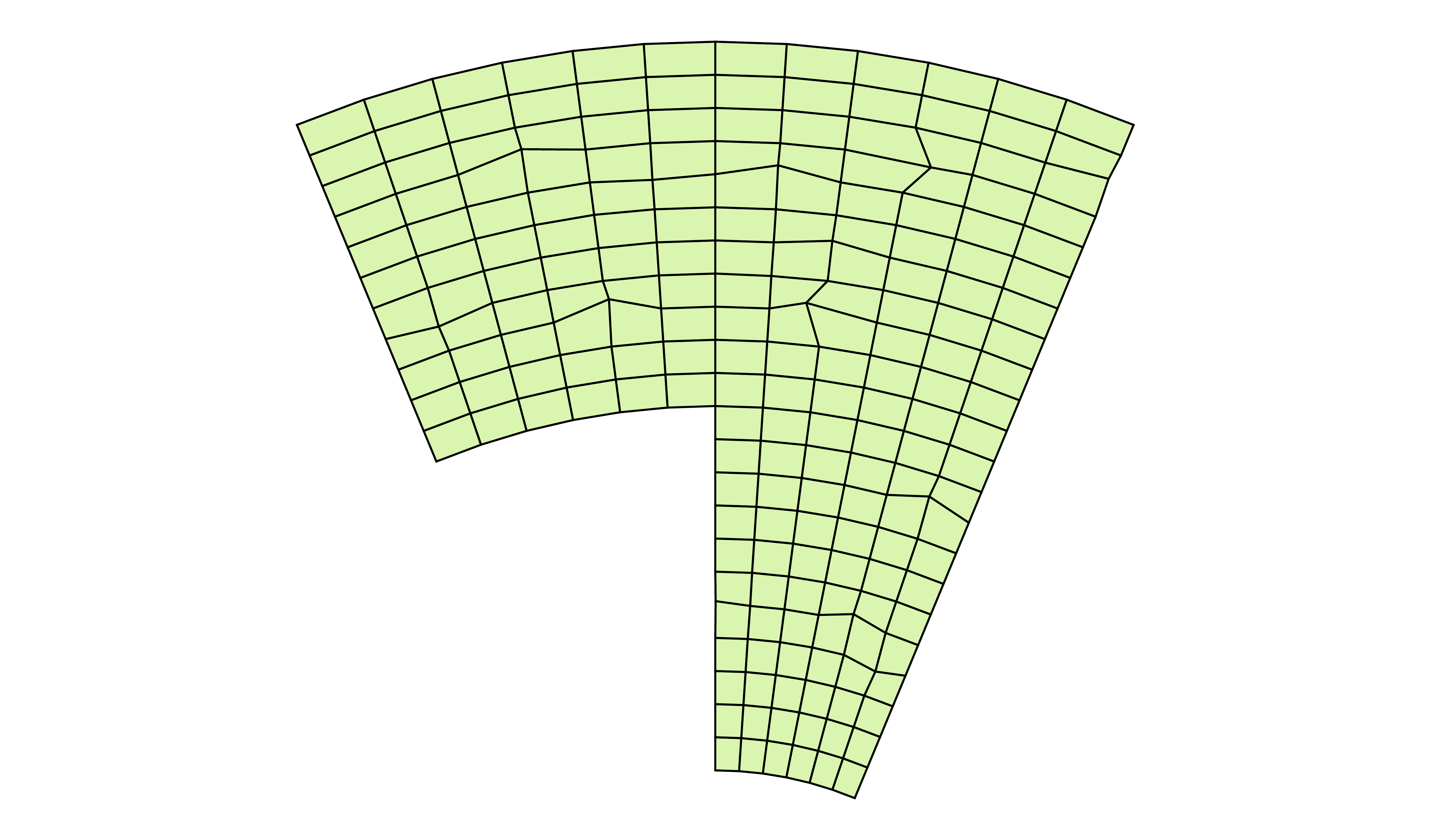}        
\caption{}
\label{dist_a}
\end{subfigure}
\begin{subfigure}{0.35\columnwidth}
\centering
\includegraphics[trim={16cm 0cm 16cm 0cm},width=0.41\columnwidth]{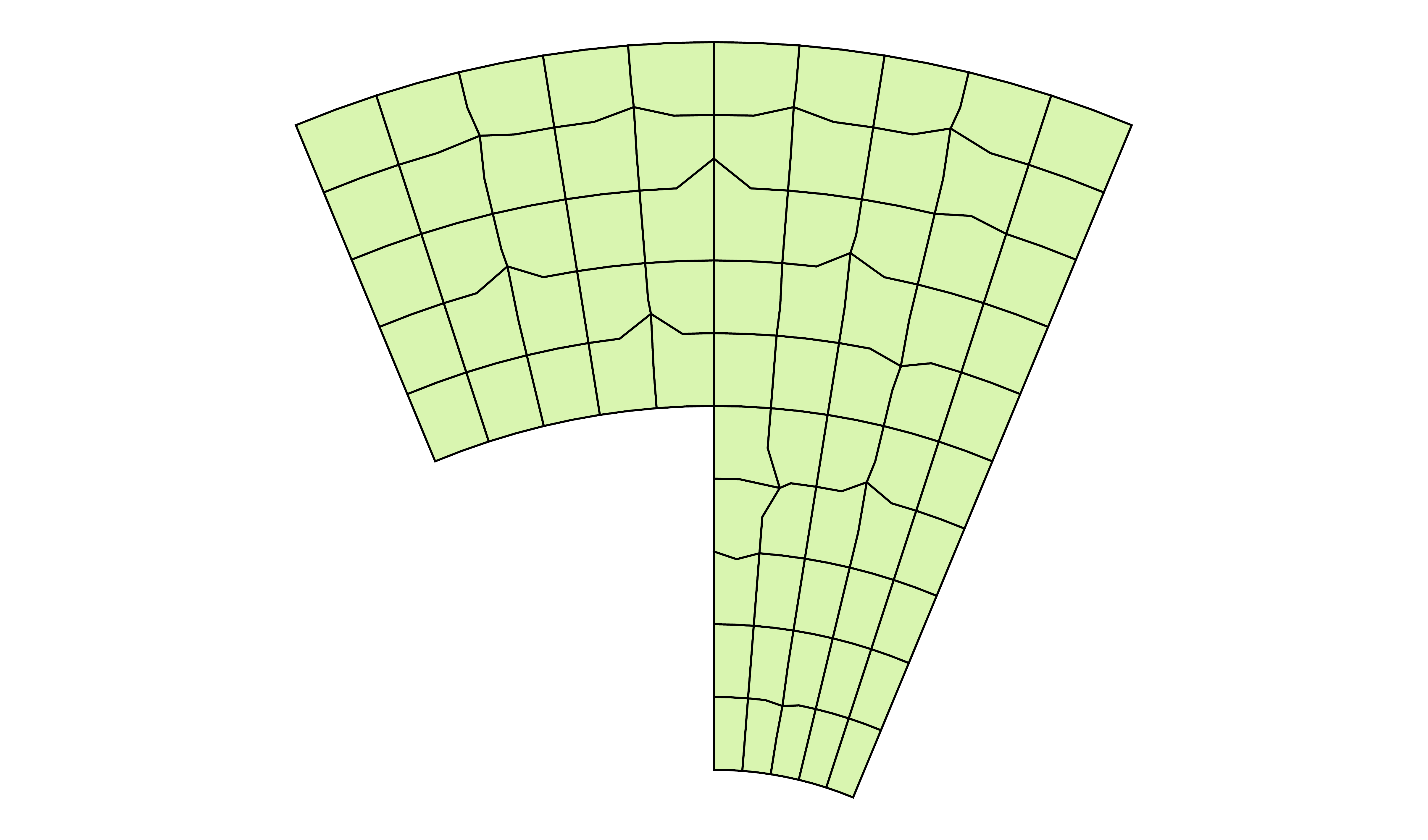}
\caption{}
\label{dist_b}
\end{subfigure}
\begin{subfigure}{0.3\columnwidth}
\centering
\includegraphics[trim={16cm 0cm 16cm 0cm},width=0.48\columnwidth]{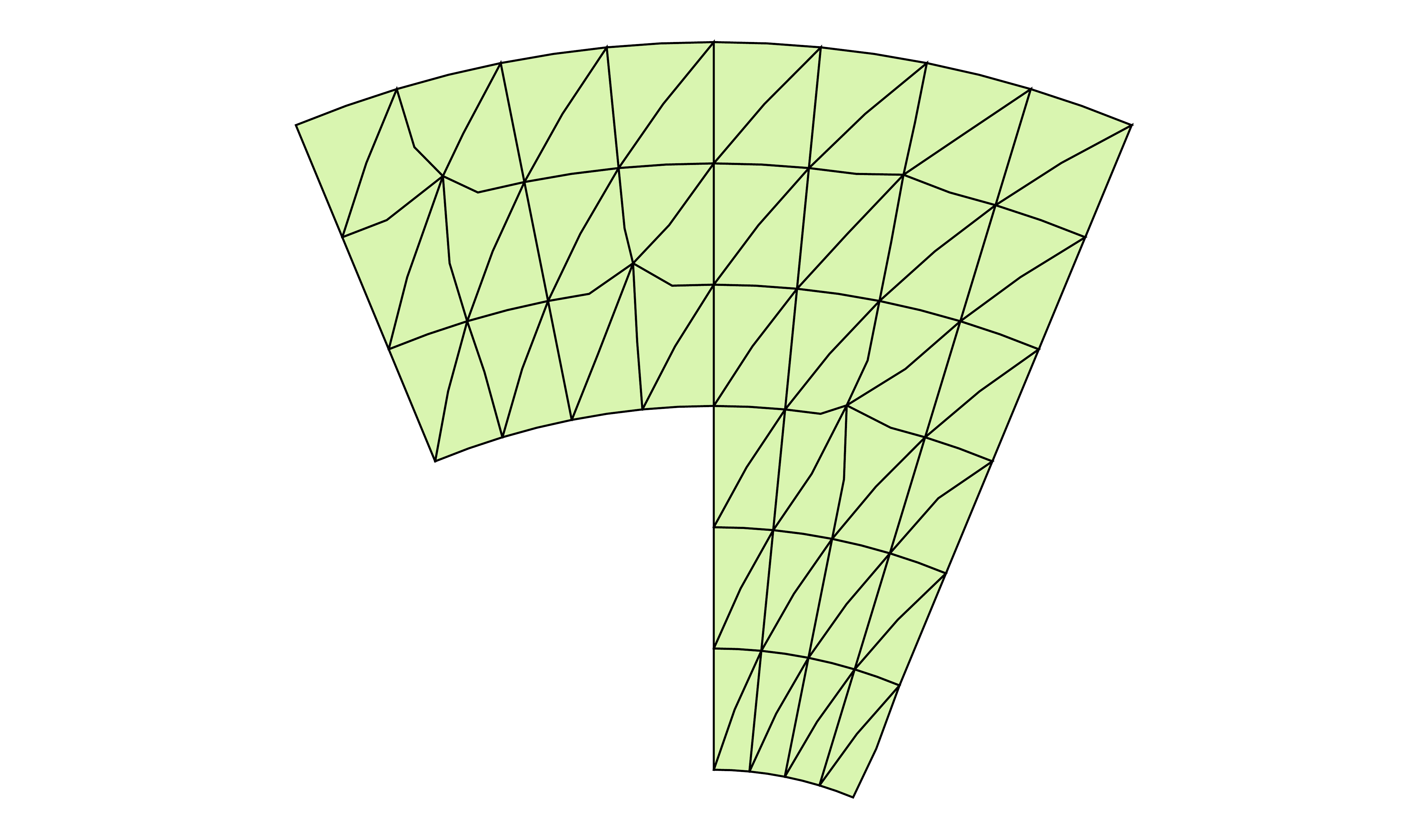}
\caption{}
\label{dist_c}
\end{subfigure}
\caption{Curved L-shaped domain discretized with (a) distorted Q4 elements, (b) distorted Q9 elements and (c) distorted T6 elements.}
\label{distorted_curvedl}
\end{figure}

 To analyze the performance of the distorted mesh we have considered the standard curved L shape problem. To obtain the eigenvalues, the curved L shape domain is discretized using a uniform and deformed mesh made of 198 Q4 and EQ4 elements, 75 Q9 and EQ12 elements and 72 T6 and ET8 elements. The distorted meshes of different elements are shown in Figs.~\ref{dist_a},~\ref{dist_b} and~\ref{dist_c}. All of the elements' results are presented in Table~\ref{tabcurlshape}. The results of a distorted mesh of all the elements are almost the same as the results of a uniform mesh. 

\vspace{5cm}
\section{Conclusions}
\label{conclusions}
In this work, a detailed performance comparison between nodal and edge elements has been presented. Various solved examples include all possible complexities like curved boundaries, non-convex domains, sharp corners, distorted meshes, and non-homogeneous domains. In every case, edge elements have shown better coarse mesh accuracy than nodal elements. We have observed in many cases that, in order to achieve the same level of accuracy, the required no. of equations with edge element is less than half of that with nodal elements.
For the non-convex domains with sharp corners, nodal elements can not predict the singular eigen value which is well predicted by all the edge elements. In addition, for such domains, nodal elements predict one additional spurious eigen value which is not present with edge elements. Also, we have observed that mesh distortion does not affect the performance of both nodal and edge elements.



%
%

\bibliographystyle{cas-model2-names}

\bibliography{library.bib}

\end{document}